\documentclass[a4paper,fleqn]{cas-sc} 

\usepackage{textcomp}
\usepackage{stackrel}
\usepackage[edges]{forest}
\usetikzlibrary{shadows.blur}
\usepackage{rotating}
\usepackage{tabularx}
\usepackage{arydshln}
     \usepackage{pdflscape}

\usepackage{rotating}
\usepackage{nicematrix}
\NiceMatrixOptions{
code-for-first-row = \color{blue} ,
code-for-last-row = \color{blue} ,
code-for-first-col = \color{blue} ,
code-for-last-col = \color{blue}
}
\usepackage{multirow}
\usepackage{array}
\usepackage{rotating,tabularx}
\usepackage{booktabs}
\usepackage{rotating}
\usepackage{color}
\usepackage{subcaption}
\usepackage{tikz}
\usetikzlibrary{matrix}
\usepackage{tabulary}
\usepackage{longtable}
\usepackage{caption}
\usepackage{booktabs,siunitx}
\usepackage{listings}
\usepackage[toc,page]{appendix}
\usepackage{float}
\usepackage{enumitem}
\usepackage[utf8]{inputenc}
\usepackage{listings,lstautogobble}
\usepackage{makecell}
\usepackage{caption}
\usepackage{subcaption}
\usepackage{tabulary}
\usepackage{algorithmic, algorithm} 
\usepackage{fancyhdr}
\usepackage{amssymb}
\usepackage{amsthm}
\usepackage{placeins}
\usepackage{amsfonts}
\usepackage[all]{xy}
\usepackage{tikz}
\usepackage{verbatim}
\usepackage{multirow}
\usepackage{psfrag}
\usepackage{pdflscape}

\usetikzlibrary{shapes,fit,positioning,shadows,arrows.meta}
\tikzset{
     auto,node distance =1 cm and 1 cm,semithick,
     color=black
}

\tikzset{parent/.style={align=center,text width=2cm,fill=green!20,rounded corners=2pt},
    child/.style={align=center,text width=2.8cm,fill=green!50,rounded corners=6pt},
    grandchild/.style={fill=pink!50,text width=2.3cm}
}
\usepackage{siunitx}
\usepackage{fancyhdr}
\usepackage{amssymb}
\usepackage{amsthm}
\usepackage{placeins}
\usepackage{amsfonts}
\usepackage[all]{xy}
\usepackage{tikz}
\usepackage{verbatim}
\usepackage{multirow}
\usepackage{psfrag}
\usepackage[authoryear,longnamesfirst]{natbib}
\hypersetup{ colorlinks, citecolor=blue, linkcolor=blue, urlcolor=Blue}
\usepackage{pdflscape}
\usepackage{multirow}
\usepackage{url}
\usepackage[utf8]{inputenc}
\usepackage{booktabs}
\usepackage{color}
\usepackage{colortbl}
\usepackage{multirow}
\usepackage{tikz}
\usepackage{booktabs}
\usepackage{rotating}
\usepackage{color}
\usepackage{tikz}
\usetikzlibrary{matrix}
\usepackage{booktabs,siunitx}
\usepackage{listings}
\usepackage[toc,page]{appendix}
\usepackage{adjustbox}
\usepackage{float}
\usepackage{enumitem}
\usepackage[utf8]{inputenc}
\usepackage{listings,lstautogobble}
\usepackage{makecell}
 \usepackage{url}
\usepackage{tabulary}
\usepackage{fancyhdr}
\usepackage{amssymb}
\usepackage{amsthm}
\usepackage{placeins}
\usepackage{amsfonts}
\usepackage[all]{xy}
\usepackage{tikz}
\usepackage{verbatim}
\usepackage{multirow}
\usepackage{psfrag}
\usepackage{pdflscape}
\usepackage{pgfplots}
\definecolor{bblue}{HTML}{4F81BD}
\definecolor{rred}{HTML}{C0504D}
\definecolor{ggreen}{HTML}{9BBB59}
\definecolor{ppurple}{HTML}{9F4C7C}
\definecolor{yellow}{HTML}{E9D66B}
\definecolor{app}{HTML}{D7191C}
\definecolor{method}{HTML}{FDAE61}
\usepackage{amsmath}
\usepackage{booktabs,siunitx}
\usepackage{listings}
\usepackage[toc,page]{appendix}
\usepackage{float}
\usepackage{enumitem}
\usepackage[utf8]{inputenc}
\usepackage{listings,lstautogobble}
\usepackage{makecell}
\newenvironment{rcases}
  {\left.\begin{aligned}}
  {\end{aligned}\right\rbrace}
\definecolor{RoyalBlue}{cmyk}{1, 0.50, 0, 0}
\usetikzlibrary{positioning, fit, calc}   
\tikzset{block/.style={draw, thick, text width=2cm , minimum height=1.3cm, align=center},   
line/.style={-latex}     
}  

\def\tsc#1{\csdef{#1}{\textsc{\lowercase{#1}}\xspace}}
\tsc{WGM}
\tsc{QE}
\tsc{EP}
\tsc{PMS}
\tsc{BEC}
\tsc{DE}
\usepackage{multicol}
\usepackage{nomencl}
\renewcommand\nompreamble{\begin{multicols}{2}}
\renewcommand\nompostamble{\end{multicols}}
\makenomenclature
\begin{document}

\shorttitle{Survey on Optimization Methods for Cyber-physical Networks}
\shortauthors{Aslani et al.}

\title [mode = title]{The State-of-the-Art Survey on Optimization Methods for Cyber-physical Networks}                      
\author[1]{Babak Aslani}[orcid=0000-0002-2734-9398]
\ead{baslani@gmu.edu}
\address[1]{Department of Systems Engineering and Operations Research, George Mason University, Fairfax, VA}

\author[1]{Shima Mohebbi}[orcid=0000-0002-1587-0506]
\ead{smohebbi@gmu.edu}
  \cormark[1]

\author[2]{Edward J. Oughton}[orcid=0000-0002-2766-008X]
\ead{eoughton@gmu.edu}
\address[2]{Department of Geography  and Geoinformation Science, George Mason University, Fairfax, VA}







\begin{abstract}
Cyber-Physical Systems (CPS) are increasingly complex and frequently integrated into modern societies via critical infrastructure systems, products, and services. Consequently, there is a need for reliable functionality of these complex systems under various scenarios, from physical failures due to aging, through to cyber attacks. Indeed, the development of effective strategies to restore disrupted infrastructure systems continues to be a major challenge. Hitherto, there have been an increasing number of papers evaluating cyber-physical infrastructures, yet a comprehensive review focusing on mathematical modeling and different optimization methods is still lacking. Thus, this review paper appraises the literature on optimization techniques for CPS facing disruption, to synthesize key findings on the current methods in this domain. A total of 108 relevant research papers are reviewed following an extensive assessment of all major scientific databases. The main mathematical modeling practices and optimization methods are identified for both deterministic and stochastic formulations, categorizing them based on the solution approach (exact, heuristic, meta-heuristic), objective function, and network size. We also perform keyword clustering and bibliographic coupling analyses to summarize the current research trends. Future research needs in terms of the scalability of optimization algorithms are discussed. Overall, there is a need to shift towards more scalable optimization solution algorithms, empowered by data-driven methods and machine learning, to provide reliable decision-support systems for decision-makers and practitioners.  

\end{abstract}


\maketitle

\section{Introduction}
Cyber-physical systems (CPS) integrate computation, communication, sensor, and network technologies equipped with feedback loops for managing interdependency among their physical and cyber components \citep{gurdur2018systematic, tu2019robustness}. Indeed, the modern networked world in which we live has made it possible to have near-ubiquitous information at hand, leading to the increased embedding of this information into CPS approaches \citep{jazdi2014cyber}. CPS have emerged by integrating Information and Communication Technology with Critical Infrastructures such as water, energy, and transportation systems, leading to increased complexity \citep{oughton_infrastructure_2018}. In specific, water infrastructure includes physical elements (pipeline, pump stations, detention basins, treatment facilities), cyber elements (SCADA systems). Similarly, transportation infrastructures entail physical (roads, railways), cyber (traffic control technologies and sensors) elements \cite{samih2019smart}. Although this novel integration enhances the efficiency and service level of CPS, it also significantly increases the connections among system components and the interdependencies between different sectors of CPS such as the cyber-physical-social interdependency between water, transportation, and cyber infrastructure systems \citep{mohebbi2020cyber}.

Nowadays, CPS are ubiquitous, with different functionalities and capabilities, often supporting critical missions that have significant economic and societal importance. The emerging CPS, such as smart cities, autonomous vehicles, and modern transportation systems are expected to be highly intelligent, electrified, and connected \citep{broo2021cyber}. Indeed, major trends in new wireless technologies (e.g., 5G/6G) focus on enabling wide-area connectivity for remote control of previously unconnected assets \citep{oughton_surveying_2022}. Thus, this trend will have significant ramifications for CPS. While there are many definitions of what constitutes critical infrastructures, the United States has a sixteen-sector definition ranging from financial services to the defence industrial base \citep{DHS}.

In recent decades, the CPS around the world experienced disastrous situations, such as Hurricane Katrina in 2005, the earthquakes in Japan in 2011, and Hurricanes Harvey, Irma, and Maria in 2017, all of which profoundly impacted the economic growth, social development, and public safety \citep{mohebbi2020cyber}. For instance, Hurricane Irma caused a widespread power outage to customers in Florida, and Hurricane Harvey triggered a significant disruption in the transportation system of Houston, Texas \citep{rahimi2022predictive}. Given the presence of such threats to CPS, decision-makers seek to have decision-making tools (both quantitative and qualitative) to protect the standard functionality of such systems in dealing with disruptions \citep{aslani2022learn}. In particular, researchers have attempted to utilize optimization methods for proactive and restoration decisions among CPS.

The development of state-of-the-art commercial solvers, such as \textit{Gurobi} and \textit{CPLEX}, and the interpretability of results could be mentioned as main drivers for the increasing trend of optimization algorithms in this application area. Even though some review papers focus on the solution approaches and interdependencies among CPS (e.g., \cite{ouyang2014review} and \cite{mohebbi2020cyber}), a systematic evaluation of mathematical models and optimization methods developed in the context of CPS is absent in the literature. In addition, the previous reviews mainly focus on interdependent networks, while a significant share of related studies investigate standalone infrastructure systems. Therefore, a comprehensive analysis of optimization methods for CPS will provide insightful directions for future research in theory and application.   

We devised a comprehensive strategy in this study to identify the main trends in the literature and pinpoint the gaps worthy of investigation in future works. In the mathematical modeling aspect, we attempted to decompose the models based on essential elements such as objective function, decision variables, and constraints to summarize a large body of information in an interpretable format. On the other hand, we separated the literature for deterministic and stochastic models to analyze the optimization-based solution algorithms for each class regarding their specific features. We also thoroughly investigated the scale of CPS networks and failure types to highlight the missing conceptual assumptions in the literature. The complementary element of our analysis is the bibliographic analysis to explore the link between methodology-based and application-based keywords, publication outlets, and how these two aspects are interconnected in practice.  

The main contribution of this paper is providing a comprehensive review of optimization methods and mathematical modeling for CPS. We explored all major databases prudently to select the most reliable pool of papers based on the scope and methodology. The result of the search process, 108 articles, has been organized and analyzed based on the formulation type, optimization solution approach, modeling elements, network size, and failure types. The review focuses on the studies developing mathematical models and applying an optimization method to solve the underlying problem. In addition, we utilize a holistic approach to relate the solution algorithm design, CPS network features, and conceptual/modeling practices to provide insights into developing practical decision support systems with enriched optimization frameworks for city-scale CPS.

The remainder of the paper is structured as follows: In section \ref{def}, we present the definition of fundamental concepts and technical terms. Next, section \ref{search} explains the search process and selection criteria to establish the pool of literature to be reviewed. In section \ref{class}, we categorize the selected studies based on multiple criteria. Additionally, section \ref{analyze} includes detailed solution methods and mathematical modeling analyses. Finally, the discussion and future research directions are provided in section \ref{discussion}. 

\section{Definition of Concepts and Assumptions} \label{def}
\subsection{Mathematical Modeling}
Mathematical modelling is the art of capturing and describing a real-world problem in the form of mathematical terms, usually in the form of equations, and then using these equations both to help understand the original problem, and also to discover new features about the problem. The modeling process can be construed as an iterative process in which real-life problems are translated into mathematical language, solved within a symbolic system, and the solutions is implemented within the real-life problem environment. Mathematical modeling aims to capture the essential characteristics of a complex real-life problem and transform them into a more abstract representation \citep{keller2017mathematical}. In this work, we focus on optimization problems where the objective function(s) and constraints for a system can be expressed by linear or non-linear functions in the form of standard mathematical models.  

\subsubsection{Deterministic Models}
In a deterministic mathematical model, it is assumed that a set of all variable states is uniquely determined by parameters in the model and by sets of previous states of these variables. In other words, the parameters and variables of the system are fully known and invariable over the optimization horizon. The general form of a deterministic model will be as follows:

\begin{equation}
\begin{cases}
 \mathop{minimize} f(x)= c_1 x_1+...+c_n x_n \\
\textbf{ subject to }:  a(x_i) \geq b_i , i=1,2,...,p\\
x_i\geq 0, j=1,2,...,n\\
c \in  C ,a \in  A ,b \in  B ,x \in  X  
\end{cases}
\end{equation}

Where the information about tuple $\langle{A|B|C}\rangle$ is fully known prior to the optimization and assumes to be unvarying over the time horizon.

\begin{itemize}
    \item \textbf{Single-objective Models}: 
A single-objective optimization problem refers to a problem where a number of decision variables $x_i$ are defined to minimize (or maximize) one desired objective function $f$ over a set of constraints $C$. 

\begin{equation}
 \mathop{minimize}_{x} \{f(x) \textit{ subject to } x \in  X   \subseteq \mathbb{R}^n \}
\end{equation}
where $X$ is the feasible set of solutions for the optimization problem.

We can define the single-objective optimization problem with both equality and inequality set of constraints. Formally, a minimization problem, with $p$ inequality constraints and $m$ equality constraints is represented by:  

\begin{equation}
\begin{cases}
\mathop{minimize} f(x) \\
\textbf{ subject to }:  g_i(x) \geq 0 , i=1,2,...,p\\
h_j(x)=0, j=1,2,...,m\\
x \in  X   \subseteq \mathbb{R}^n
\end{cases}
\end{equation}

\item \textbf{Multi-objective Models}: 
 Optimization problems in practice can include several conflicting objectives to minimize or maximize simultaneously. A general continuous multi-objective problem aims to find $n$  decision variables $x \in \mathbb{R}^n$ that simultaneously minimize (or maximize) $k$ objective functions such as $f_r: \mathbb{R}^n \longmapsto \mathbb{R}, r= 1,…,k$ \citep{kadzinski2017evaluation}. Similar to single-objective models, the decision variables and objectives are subject to a set of bounds and constraints. A feasible solution to the Multi-Objective Optimization problem satisfies all the bounds for a decision variable, together with the $p$ and $m$ inequalities and equality constraints. The standard form of a multi-objective optimization problem will be as follows:
 \begin{equation}
\begin{cases}
\mathop{minimize} f(x,p)  \equiv ( f_{1}(x),f_{2}(x),...f_{k}(x))\\
\textbf{ subject to }:  h_{i}(x)=0, i=1,2,...,m\\
  g_{j}(x) \geq 0, i=1,2,...,p\\
  x_l \in  [x_l^{L},x_l^{U} ]  , l=1,2,...,n
\end{cases}
\end{equation}
\end{itemize}

\subsubsection{Stochastic Models}
The underlying assumption in stochastic models is the presence of uncertainty in the parameters and decision variables of the model. In other words, the variable and parameter states are not unique values but rather are presented in the form of probability distributions. Let $X$ be the domain of all feasible decisions and $x$ a specific decision. The optimization goal is to search over $X$ to find a decision that minimizes (or maximizes) an $F$ function. Let $\xi$ denote random information available only after the decision is made. In the general form of stochastic modeling, we define the objective function as a random cost function such as $F(x,\xi)$ and constraints can also be written as random function. Since $F(x,\xi)$ cannot be optimized directly, the focus will be on the expected value, $ \mathbb{E}[F(x,\xi)]$. The general form of a stochastic optimization problem can be presented as follows \citep{hannah2015stochastic}:

\begin{equation}
    \zeta^*= \mathop{minimize}_{x \in X} \{f(x)= \mathbb{E}[F(x,\xi)] \}
\end{equation}

Where the optimal solution will be the set of $S^*= \{x \in X : f(x) = \zeta^*\}$.

\begin{itemize}
    \item \textbf{Two-stage stochastic models}: Two-stage stochastic mathematical model is the most common formulation in the field of stochastic programming. The basic idea of this modeling approach is that optimal decisions should be based on data available at the time the decisions are made and cannot depend on future observations. A general two-stage stochastic programming problem is given by \citep{shapiro2007tutorial}:

\begin{equation}
{\displaystyle {\begin{array}{llr}\min \limits _{x\in \mathbb {R} ^{n}}&g(x)=c^{T}x+E_{\xi }[Q(x,\xi )]&\\{\text{subject to}}&Ax=b&\\&x\geq 0&\end{array}}}
\end{equation}

where ${\displaystyle Q(x,\xi )}$ is the optimal value of the second-stage problem given below:
\begin{equation}
{\displaystyle \min _{y}\{q(y,\xi )\,|\,T(\xi )x+W(\xi )y=h(\xi )\}.}
\end{equation}


In such formulation, ${\displaystyle x\in \mathbb {R} ^{n}}$ is the first-stage decision variable vector, ${\displaystyle y\in \mathbb {R} ^{m}}$ is the second-stage decision variable vector, and ${\displaystyle \xi (q,T,W,h)}$ contains the data of the second-stage problem. In the first stage, we have to make a "here-and-now" decision ${\displaystyle x}$ before the realization of the uncertain data ${\displaystyle \xi }$, viewed as a random vector, is known. In the second stage, after a realization of ${\displaystyle \xi }$ becomes available, we optimize the decisions by solving a new optimization problem updated with the new information. The general case of a two-stage problem is linear with linear terms for objective function and the constraint. However, conceptually this is not essential, and integer or nonlinear assumptions could also be incorporated into two-stage models \citep{pichler2016nonlinear}.

\item \textbf{Multi-stage stochastic models}: This stochastic formulation has been introduced as a generalization of two-stage models to capture the sequential realization of uncertainties and labelled as \textit{multistage stochastic programming}. In this setting, the time horizon is discretized into “stages” where each stage has the realizations of uncertainties at the current stage. A wide range of complex and dynamic optimization problems are effectively modeled as multi-stage stochastic programs, possessing key decision variables in one or more of the stages. In the multi-stage setting, the uncertainty of the problem is revealed sequentially in $T$ stages. The random variable $\xi$ is split into $T-1$ chunks, $\xi=(\xi_1,...,\xi_{T-1})$, and the problem at hand is to decide at each stage $t=1,...,T$ what is the optimal action, $x_t(\xi_{[1,t-1]})$, given the previous observations $\xi_{[1,t-1]}:=(\xi_1,...,\xi_{t-1})$. The global variable of this problem can be written as follows \citep{bakker2020structuring}:

\begin{equation}
    x(\xi)=(x_1,x_2(\xi_1), \dots, x_T(\xi_{[1,T-1]})) \in \mathbb{R}^{n_1} \times \dots \times \mathbb{R}^{n_T}
\end{equation}

where $(n_1,...,n_T)$ are the size of the decision variable at each stage and $n=\sum \limits_{t=1}^T n_t$ is the total size of the problem.

\item \textbf{Markov Decision Process (MDP)}: MDP is a mathematical framework for sequential decision making under uncertainty
that has informed decision making in a discrete-time stochastic environment. This branch of stochastic programming has been adopted in a variety of application areas including inventory control, scheduling, and medicine \citep{steimle2021multi}. 

A Markov decision process is formulated in the form of a tuple such as ${\displaystyle \langle S,A,P_{a},R_{a}\rangle }$, where:
\begin{itemize}
    \item  ${\displaystyle S}$ is a set of states or the state space.
    \item ${\displaystyle A}$ is a set of actions or the action space.
    \item ${\displaystyle P_{a}(s,s')=\Pr(s_{t+1}=s'\mid s_{t}=s,a_{t}=a)}$ is the probability that action ${\displaystyle a}$ in state ${\displaystyle s}$ at time ${\displaystyle t}$ will lead to state ${\displaystyle s'}$ at time ${\displaystyle t+1}$.
    \item ${\displaystyle R_{a}(s,s')}$ is the immediate reward (or expected immediate reward) received after transitioning from state ${\displaystyle s}$ to state ${\displaystyle s'}$, due to action ${\displaystyle a}$.
\end{itemize}
 
A policy function ${\displaystyle \pi }$  is a probabilistic mapping from state space (${\displaystyle S}$) to action space (${\displaystyle A}$).

At each time step, the process is in some state ${\displaystyle s}$, and the decision-maker may choose any action ${\displaystyle a}$ that is available in state ${\displaystyle s}$. The process responds at the next time step by randomly moving into a new state ${\displaystyle s'}$, and giving the decision-maker a corresponding reward ${\displaystyle R_{a}(s,s')}$. The probability that the process moves into its new state ${\displaystyle s'}$ is influenced by the chosen action; however, it is independent of all previous states and actions to satisfy the Markov property.

The goal of a Markov decision process is to find an optimal \textit{policy} for the decision-maker: a function ${\displaystyle \pi }$ that specifies the action ${\displaystyle \pi (s)}$ that the decision-maker will choose when in state ${\displaystyle s}$. The objective of a general MDP is to choose a policy ${\displaystyle \pi }$  that will maximize some cumulative function of the random rewards over a infinite horizon:
\begin{equation}
{\max (Z)=\displaystyle E\left[\sum _{t=0}^{\infty }{\gamma ^{t}R_{a_{t}}(s_{t},s_{t+1})}\right]}
\end{equation}

Where ${\displaystyle a_{t}=\pi (s_{t})}$ is the selected actions given by the policy, and the expectation is taken over ${\displaystyle s_{t+1}\sim P_{a_{t}}(s_{t},s_{t+1})}$. In this equation, $\gamma$ is a discount factor satisfying $ 0\leq \gamma \leq  1$, which is usually close to 1.
\end{itemize}

\subsection{Optimization Solution Algorithms}
Global optimization algorithms can be broadly divided into three groups including \textit{exact methods} (e.g., branch-and-bound, cutting planes, and decomposition) and \textit{heuristic methods}, and \textit{evolutionary methods} (e.g., Genetic Algorithm (GA)). The main difference in these classes originates in the search process and the optimality guarantee. In other words, while in the exact methods can provide solution with guaranteed optimality, there is no certainty about the optimal status of the output for heuristic and evolutionary methods.  

\begin{itemize}
    \item \textbf{Exact methods}: Exact methods can comb the solution space of any optimization problem to find the optimal solution among all candidate feasible solutions. This guaranteed optimality is the primary rationale for implementing these approaches for different classes of optimization problems. An exact optimization method is essentially a solution method of choice that can solve an optimization problem with an effort that grows polynomially with the problem size. The most common classes are branch-and-bound, cutting plane, and decomposition methods.
    
    \item \textbf{Heuristic methods}
    Approximation algorithms are a diverse family of solution methods for optimization problems that return an approximate solution without any promise of optimality. Specifically, heuristics solution methods are problem-specific approximate algorithms attempting to exploit the underlying rules of the problem at hand to obtain a high-quality solution. Even though heuristics do not guarantee to find an optimal solution, they are simple and much faster than exact approaches. There are two different types of heuristics: \textit{Construction heuristics}, which construct one solution from scratch by performing iterative construction steps, and \textit{Improvement heuristics}, which start with a complete solution and iteratively apply problem-specific search operators to search the solution space. Among the well-known heuristics, nearest neighbor heuristic, cheapest insertion heuristic, and \textit{k}-opt heuristic are noteworthy that are developed for Travelling Salesman Problem (TSP) problem but are applicable to other optimization problems as well \citep{rothlauf2011optimization}
        
    \item \textbf{Evolutionary methods}
    Evolutionary algorithms (or interchangeably Meta-heuristics) are computationally efficient methods for solving complex optimization problems approximately, particularly for the Combinatorial Optimization Problems class of optimization problems with discrete decision variables and finite search space. Since these methods can provide near-optimal solutions in reasonable computational time for NP-Hard class of optimization problems, they are often considered a good substitution for exact solution algorithms \citep{Karimi2022}. From a technical perspective, evolutionary methods are approximate optimization methods constructed based on the interaction between local search procedures and higher-level search strategies to generate an iterative search process able to escape from local optimality and perform a comprehensive exploration of the search space \citep{gendreau2010handbook}. Typical examples of evolutionary algorithms are GA, Simulated Annealing (SA), Tabu Search (TS), and Greedy Randomized Adaptive Search Procedure (GRASP).

\end{itemize}

\subsection{Machine Learning}
Machine learning (ML), a sub-field of artificial intelligence, is referred to a set of algorithmic and statistical methods attempting to \textit{learn} from data to improve the performance of computational methods in solving challenging classes of problems \citep{bishop2006pattern}. These intelligent methods autonomously improve their learning quality over time, using extracted knowledge from data and information through observations and real-world interactions. Assuming we have a set of input variables $x_i \in  X$ that are used to determine an output variable $y$,  The goal of ML is to quantify the mapping function $f$ between the input variables and the output variable in such as $Y=f(\{x_i \in  X\})$.

ML algorithms can be classified into three main groups \cite{bishop2006pattern}:

\begin{itemize}
    \item \textbf{Supervised Learning Algorithms}: In this class of ML methods, the values of input variables $x_i \in  X$ and the corresponding values of the output variables $y$ are known a priori. A supervised learning algorithm aims to automatically discover the relationship between input variables and output to predict the outcome $y^\prime$ for new input variables $ X^\prime$. The supervised learning algorithms can be classified into two sub-classes of classification and regression algorithms based on their scope. Classical supervised learning algorithms include Linear Regression, Support Vector Machine, Naive Bayes , Gradient Boosting, Decision Tree, Random Forest, and k-Nearest Neighbor.
    
    \item \textbf{Unsupervised Learning Algorithms}: The ML methods in this class are adopted when the values of input variables  $x_i \in  X$ are known while there are no associated values for the output variables. The learning goal of these techniques is to detect and utilize the hidden patterns in the input data. Classical unsupervised learning algorithms include k-means clustering, Self-Organizing Map, and Principal Component Analysis.
    
    \item \textbf{Reinforcement Learning (RL) algorithms}: In these more recent ML algorithms, an agent learns from interactions with the surrounding environment to take actions that maximize the pre-defined reward iteratively. At each iteration, the agent automatically decides the optimal action within a specific context to optimize its performance based on a reward feedback mechanism. RL algorithms include Q-Learning, Learning Automata, opposition-based RL, and Deep Reinforcement Learning.
\end{itemize}

\section{Search process and selection criteria} \label{search}

To find the most reliable and diverse pool of papers, well-known scientific databases, including Elsevier, IEEE Explore, Springer, Wiley, and INFORMS, have been carefully searched. To this end, we tested a series of keywords and refined the keyword combination based on the quality of search results (in terms of relevance and diversity). Finally, the search process is conducted using the following search rule:

\textbf{\{(((“Cyber-infrastructure” OR “Cyber-physical” ) AND ("Infrastructure") AND (“Optimization” OR “Optimize” ) OR (“Restoration”)))\}} \footnote{ We also considered the British spelling of \textit{optimise} and \textit{optimisation} terms.  }

In the next step, we selected the literature review papers based on the following key points:

\begin{itemize}
    \item The paper should define a clear \textbf{mathematical model} (decision variables, objective function(s), and constraints).
    \item The paper should develop an \textbf{optimization-based solution algorithm}, including pure optimization algorithms, hybrid game theory and optimization models, and simulation-optimization methods.  
    \item The paper should focus on \textbf{cyber-physical infrastructure systems} dealing with \textbf{disruptive incidents}. However, the scope of the paper can be a combination of network design (high-level decision-making), preventive maintenance (proactive efforts), and restoration (reactive actions).    
    \item The paper should provide information about the \textbf{network size} of the case study or random instances to be analyzed from the scalability viewpoint. 
\end{itemize}

In this review, we divided the CPS networks studied in the literature in terms of size into three categories of small, medium-sized, and large-scale networks. The small class is referred to a network with less than 100 overall components (node and edges). We also denote any network with more than 100 and less than 400 components as medium-scale. Finally, without loss of generality, we consider a CPS network with more than 400 overall components (node and edges) as a large-scale network.

\subsection{Final Review Pool}
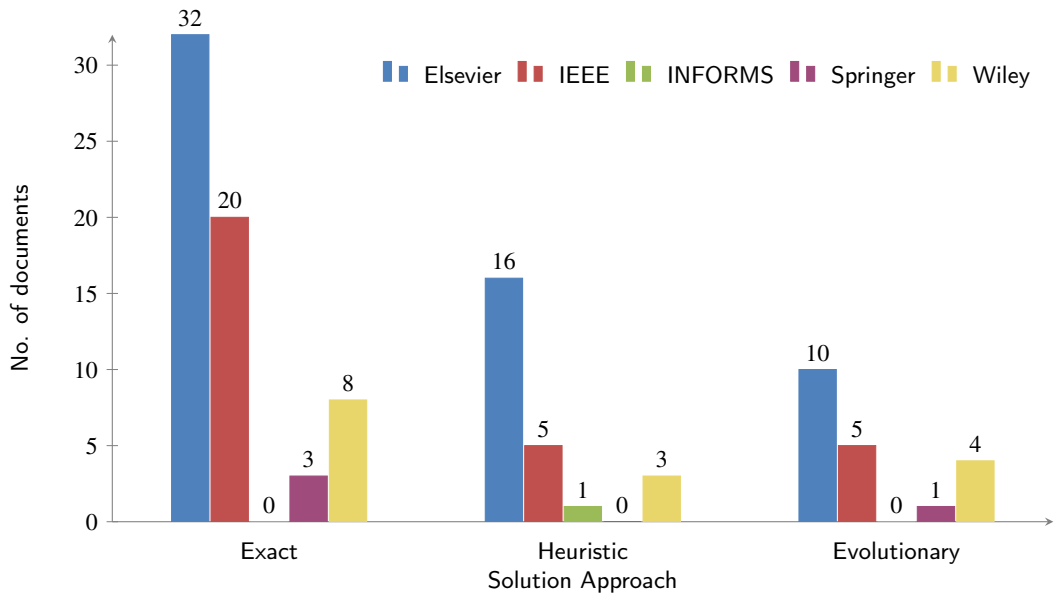
\begin{figure}[htb]
\centering
\begin{tikzpicture}
    \begin{axis}[
    grid = none,
      axis lines=left,
        width  = 0.85*\textwidth,
        height = 8cm,
        major x tick style = transparent,
        ybar=2*\pgflinewidth,
        bar width=14pt,
        nodes near coords,
        every node near coord/.append style={text=black},
        ylabel = {No. of documents},
        	y label style={at={(0,0.5)}},
         xlabel = {Solution Approach},
        symbolic x coords={Exact,Heuristic,Evolutionary},
        xtick = data,
        scaled y ticks = false,
        enlarge x limits=0.25,
        ymin=0,
        legend cell align=left,
        legend style={
        draw=none,
                at={(0.99,0.87)},
                anchor=south east,
                column sep=1ex,
                 legend columns = -1
        }
    ]
        \addplot[style={bblue,fill=bblue,mark=none}]
            coordinates {(Exact, 32) (Heuristic,16) (Evolutionary,10)};

        \addplot[style={rred,fill=rred,mark=none}]
             coordinates {(Exact,20) (Heuristic,5) (Evolutionary,5)};

        \addplot[style={ggreen,fill=ggreen,mark=none}]
             coordinates {(Exact,0) (Heuristic,1) (Evolutionary,0)};

        \addplot[style={ppurple,fill=ppurple,mark=none}]
             coordinates {(Exact,3) (Heuristic,0) (Evolutionary,1)};
             
        \addplot[style={yellow,fill=yellow,mark=none}]
             coordinates {(Exact,8) (Heuristic,3) (Evolutionary,4)};
        \legend{Elsevier,IEEE,INFORMS,Springer, Wiley}
    \end{axis}
\end{tikzpicture}
    \caption{Summary of the publisher sources for the final selected documents}
        \label{data}
\end{figure}

After scrutinizing the search results based on the outlined criteria, 108 papers are selected for the pool of literature to be reviewed. The final set of papers formulated a mathematical model for an infrastructure network dealing with disruption. In addition, all studies designed an optimization solution algorithm and applied their framework to real-life case studies or artificially generated networks. Fig. \ref{data} shows the distribution of the selected papers in each solution approach class and the database sources. 

\subsection{Mission Areas}
\begin{figure}[htbp]
	\centering
	\includegraphics[scale=0.55]{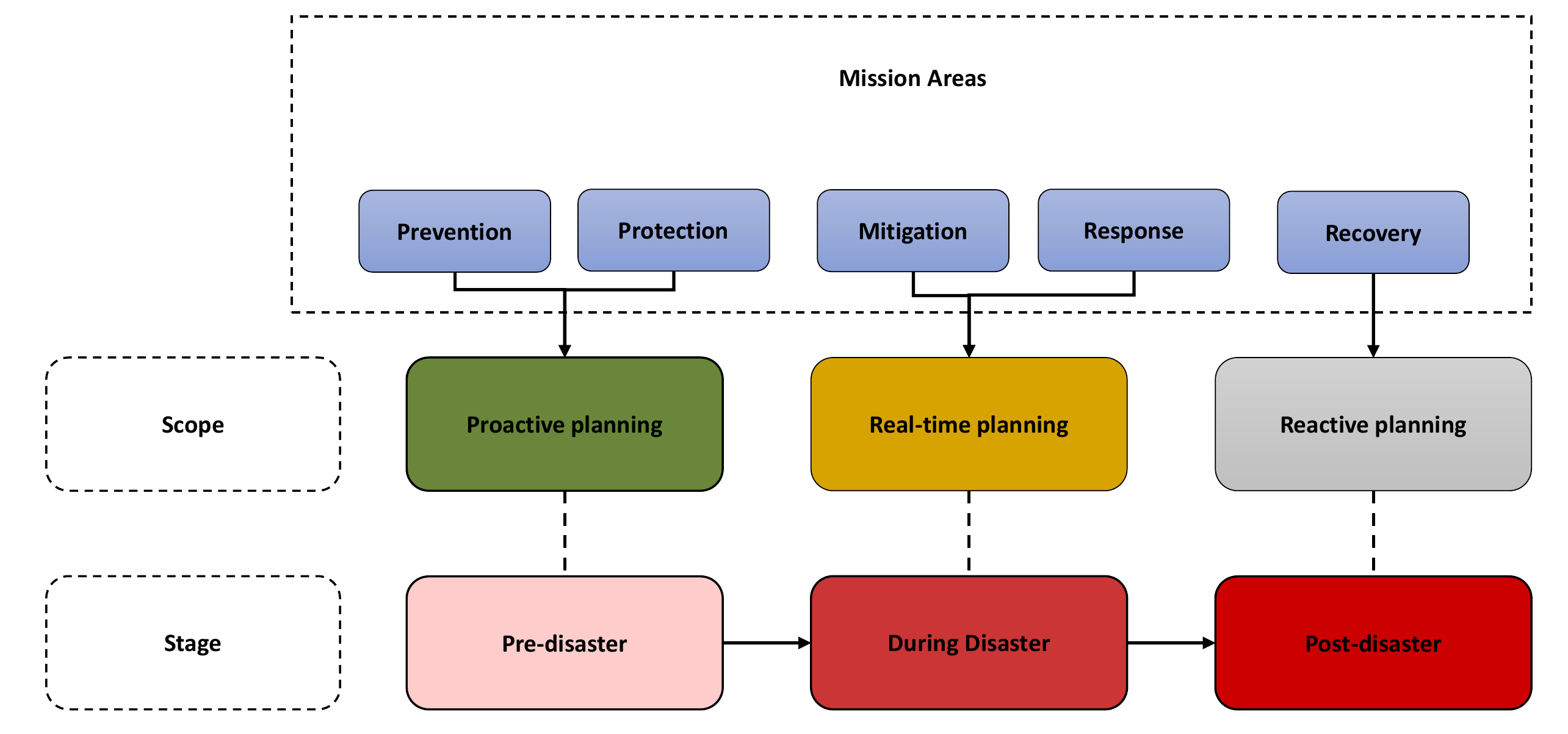}
	\centering
	\caption{Mission areas mapping for the optimization frameworks in the context o CPSs}
	\label{mission}
\end{figure}

The term \textit{mission area} is defined by the United States Department of Homeland Security and Federal Emergency Management Agency. It is one of five different classes of corrective actions capable of enhancing the preparedness of infrastructure systems in dealing with threats and hazards that pose the greatest risk to the security of the United States \citep{directive2011presidential}. In this context, the concept of mathematical modeling and optimization for CPS can cover these distinct mission areas in different phases as depicted in Fig. \ref{mission}:
\begin{itemize}
    \item \textbf{Proactive  planning}: When the optimization problem is defined in a pre-disaster situation and the actions are planned to be taken before a disruptive incident. This category covers the \textit{prevention/protection} mission areas defined in the DHS guideline. 
      \item \textbf{Real-time planning}: When the optimization problem is defined during a disaster and the actions are planned and taken during or immediately after a disruption. This category covers the \textit{mitigation/response} mission areas defined in the DHS guideline.  
        \item \textbf{Reactive planning}: When the optimization problem is defined in a post-disaster situation and the mitigating actions are planned to be taken completely after the impacts of an extreme event. This category covers the \textit{recovery} mission areas defined in the DHS guideline.  
\end{itemize}

The focus of this literature review is on the reactive planning scope in the post-disaster phase. Therefore, all the selected documents defined a problem and implemented optimization techniques to restore affected CPSs to a normal service level following a disruption.   

\section{Classifications of selected studies} \label{class}
\begin{figure}[htbp]
\centering
\begin{center}
\begin{tikzpicture}[scale=1.3]
    \begin{axis}[
    grid = none,
    axis lines=left,
    enlarge x limits=0.03,
    xtick=data,
    x tick label style={/pgf/number format/1000 sep=,rotate=90,anchor=east, font=\footnotesize},
    y tick label style={/pgf/number format/1000 sep=,anchor=east, font=\footnotesize},
    xmin = 2007, xmax = 2022,
    ymin = 0, ymax = 21,
    height=6 cm,
	width=12 cm,
	xlabel={Year},
	x label style={at={(axis description cs:0.5,-0.1)},anchor=north, font=\footnotesize},
	ylabel style={font=\footnotesize},
	ylabel={No. of documents}, 
	y label style={at={(0.05,0.5)}},
	legend style={draw=none,at={(0.22,0.97)},anchor=north, font=\footnotesize},
]
\addplot [-,smooth,color=CadetBlue,mark=*]coordinates{
(2007,2)
(2008,1)
(2009,4)
(2010,2)
(2011,1)
(2012,1)
(2013,5)
(2014,4)
(2015,9)
(2016,5)
(2017,9)
(2018,14)
(2019,17)
(2020,12)
(2021,19)
(2022,3)
};
 \addlegendentry{Overall}
 \addplot [-,dashed,color=Peach,mark=square*, mark options={solid}]coordinates{
(2007,1)
(2008,0)
(2009,3)
(2010,2)
(2011,1)
(2012,1)
(2013,4)
(2014,4)
(2015,5)
(2016,4)
(2017,8)
(2018,9)
(2019,11)
(2020,8)
(2021,15)
(2022,2)
};
 \addlegendentry{Deterministic formulation}

 \addplot [-,dashed,color=Tan,mark=triangle*, mark options={solid}]coordinates{
(2007,1)
(2008,1)
(2009,1)
(2010,0)
(2011,0)
(2012,0)
(2013,1)
(2014,0)
(2015,4)
(2016,1)
(2017,1)
(2018,5)
(2019,6)
(2020,4)
(2021,4)
(2022,1)
};
 \addlegendentry{Stochastic formulation}
\end{axis}
\end{tikzpicture}
\caption{Number of papers per year and type of mathematical formulation}
\label{time}
\end{center}
\end{figure}
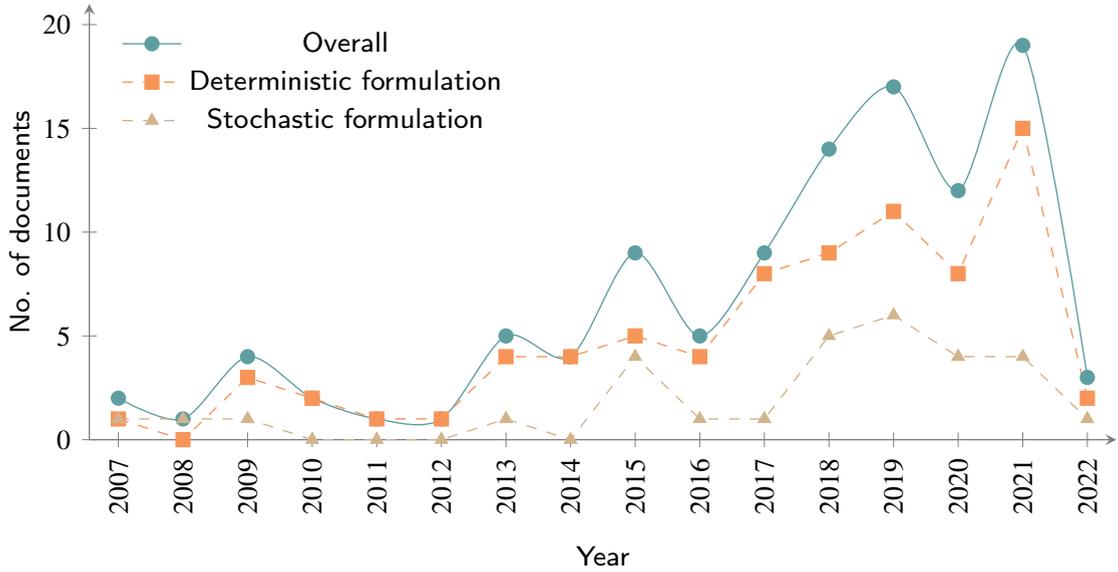

In the first classification step, we separated the studies based on the modeling assumptions into deterministic and stochastic formulations. Fig. \ref{time} shows the number of papers for each year along with the individual numbers for each modeling approach during this time horizon. We then divided the studies in each formulation class based on the proposed optimization solution approaches (i.e., exact, heuristic, and evolutionary).

\subsection{Deterministic Formulation}
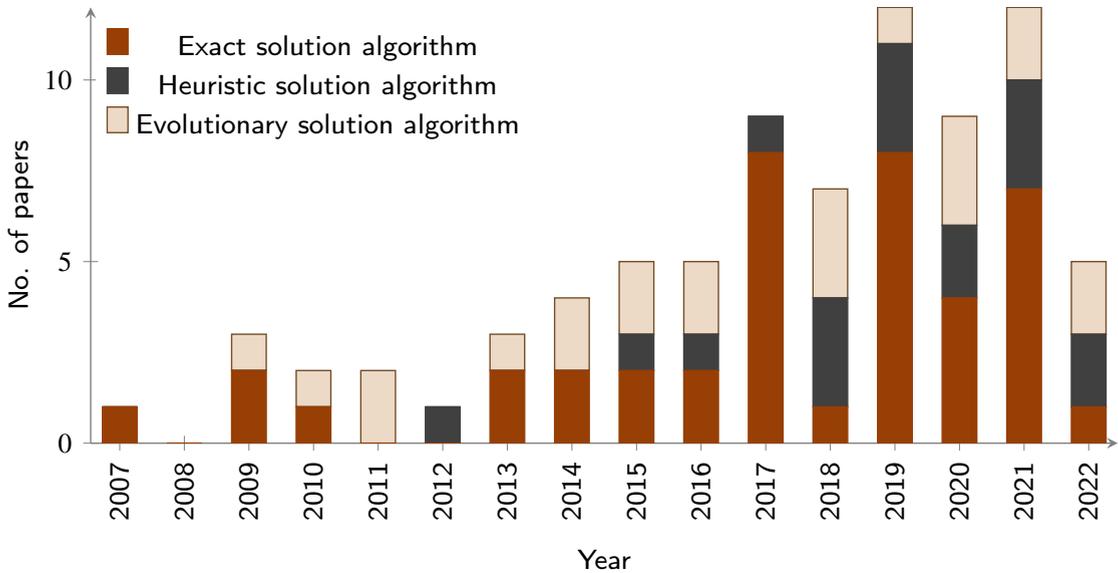
\begin{figure}[htbp]
    \centering

\begin{center}
\begin{tikzpicture}[scale=1.3]
\begin{axis}[
    ybar stacked,
    axis lines=left,
    enlarge x limits=0.03,
    ymin=0,
    x tick label style={/pgf/number format/1000 sep=,rotate=90,anchor=east, font=\footnotesize},
    y tick label style={/pgf/number format/1000 sep=,anchor=east, font=\footnotesize},
    legend style={draw=none,at={(0.22,0.97)},anchor=north, font=\footnotesize},
    height=6 cm,
	width=12 cm,
    x label style={at={(axis description cs:0.5,-0.10)},anchor=north, font=\footnotesize},
	xlabel={Year},
	ylabel style={font=\footnotesize},
    ylabel={No. of papers},
    y label style={at={(0.05,0.5)}},
    symbolic x coords={2007, 2008, 2009, 2010, 2011,2012,2013,2014,2015,2016,2017,2018,2019,2020,2021,2022
		},
    xtick=data,
    ]
\addplot+[style={RawSienna},ybar] plot coordinates {

(2007,1)
(2008,0)
(2009,2)
(2010,1)
(2011,0)
(2012,0)
(2013,2)
(2014,2)
(2015,2)
(2016,2)
(2017,8)
(2018,1)
(2019,8)
(2020,4)
(2021,7)
(2022,1)
};
\addlegendentry{Exact solution algorithm}
\addplot+[style={darkgray},ybar] plot coordinates {
(2007,0)
(2008,0)
(2009,0)
(2010,0)
(2011,0)
(2012,1)
(2013,0)
(2014,0)
(2015,1)
(2016,1)
(2017,1)
(2018,3)
(2019,3)
(2020,2)
(2021,3)
(2022,2)
};
\addlegendentry{Heuristic solution algorithm}
\addplot+[style={},ybar] plot coordinates {
(2007,0)
(2008,0)
(2009,1)
(2010,1)
(2011,2)
(2012,0)
(2013,1)
(2014,2)
(2015,2)
(2016,2)
(2017,0)
(2018,3)
(2019,1)
(2020,3)
(2021,2)
(2022,2)
};
\addlegendentry{Evolutionary solution algorithm}
\end{axis}
\end{tikzpicture}
\caption{ Number of papers per year and per each type of solution method for deterministic mathematical model.}
\label{det}
\end{center}
\end{figure}

In this section, we focus on the studies that developed deterministic mathematical models and categorize them based on the adopted solution approach. We present a comprehensive overview of the solution algorithms in each sub-section and some representative studies. Fig. \ref{det} shows the frequency of each solution approach for the deterministic formulation based on the publication date of papers. 

\subsubsection{Exact Solution Approach}
Exact methods can thoroughly explore the solution space to find the optimal solution for any optimization problem. Therefore, a significant share of selected studies adopted these solution algorithms to guarantee the optimality of the proposed solution for decision-makers in CPS networks facing disruption. From a mathematical modeling perspective, the studies developed conventional mathematical models, as well as bi-level, and tri-level optimization models. 

Mathematical models for exact solution methods are present in the selected studies in different flavors, including Mixed Integer Programming (MIP), Integer Programming (IP), and Linear Programming (LP). Among linear models, \cite{lee2007restoration} is the first work that proposed a mathematical model to minimize the total cost of flow and slack through time for a disrupted infrastructure system restoration problem. This model established the foundation of the subsequent variation of models in other studies. For instance, \cite{cavdaroglu2013integrating} provided a modified model for the extended problem of integrated restoration and planning of interdependent infrastructures to track the recovery trajectory over the restoration horizon. Likewise, the models presented in \cite{sharkey2015interdependent} aimed to find the optimal restoration strategy and assignment of available work crews was structurally similar to the mentioned studies. Among other similar works, \cite{gonzalez2016interdependent} developed a model for the Interdependent Network Design Problem to minimize the restoration cost of a partially disrupted system from hypothetical earthquakes in Shelby County, TN. As another example, \cite{baidya2017effective} modeled the restoration planning of physically interconnected infrastructure networks and solved the model optimally to prioritize the restoration activities of disrupted components while minimizing the number of required actions for complete recovery. The detailed list of linear models for each paper is presented in Table \ref{exact}.

However, linear models mainly fail to capture the real-life relationship among different variables appropriately. As a result, nonlinear models, while complex and computationally challenging, are valuable tools for describing more realistic assumptions. The term nonlinear is assigned to a model where the objective function or constraints include a form of nonlinearity. Among notable nonlinear mathematical models for restoring CPS, \cite{ouyang2017mathematical} proposed a tri-level model with nonlinear terms in both objective function and constraints for a worst-case analysis for the resilience of interdependent CIs under Spatially Localized Attacks. In a recent work, \cite{li2021resilience} developed a mixed-integer nonlinear program  for the recovery of an interconnected system of electricity-water-gas with functional interdependencies. In this model, the operational constraints of networks, posed by physical laws such as hydraulic characteristics of pumps, are the source of nonlinearity.

A particular case of mathematical models known as constrained quadratic is the other non-conventional modeling approach in the literature. In this category, \cite{mohagheghi2014reinforcement} is the only study that formulated the optimization problem as a model with a constrained quadratic term in the objective function. They provided a mathematical framework for power grids during natural disaster events to enhance resilience against natural disasters.

Bi-level models are also adopted in the literature to capture even more complex problems arising in CPS. Bi-level models are special mathematical models in which one problem is nested within another. In the optimization literature, the outer problem is referred to as the upper-level, and the inner problem is commonly referred to as the lower-level optimization task. These specific mathematical models can be found in a limited number of studies. For example, \cite{li2015bilevel} formulated a bi-level model for analyzing the impact of coordinated cyber-physical attacks, aiming at identifying the most damaging and undetectable physical attacks constrained by the budget limitation of the attacker. The other studies are summarized in the relevant sub-section in Table \ref{exact}. 

Unlike most studies in the literature which focused on restoring a disrupted system, tri-level MIP models were formulated to approach the protection-interdiction-restoration problem to optimally balance the vulnerability and recoverability of CPS before and after a disruption. Among the few studies in this area, \cite{ghorbani2021decomposition} proposed a tri-level model to improve the resilience of an infrastructure system against disruption by addressing the concurrent vulnerability reduction and recoverability enhancement actions in a Defender-Attacker-Defender model. The other studies that adopted a tri-level model are summarized in the relevant sub-section in Table \ref{exact}.

Regarding the optimization strategy, there are two main approaches in the literature (see Table \ref{exact}). The analysis revealed that the most common approach is implementing the developed mathematical models in commercial solvers (such as \textit{CPLEX}, \textit{Gurobi}, and \textit{AMPL}) and reporting the objective function value and computational time. We refer to this solution strategy as \textit{Direct Optimization} in all following references. The analysis of studies adopting this approach shed light on the technical limitation of direct optimization. This optimization option seems to apply to more miniature complex models and small-scale networks. 

On the other hand, large-scale optimization models are complex. Because of their high-dimensional decision space, this class of problems poses a significant challenge to conventional optimization methods. When dealing with large-scale optimization problems, one viable idea is to decompose the original large-scale problem into a series of smaller and simpler sub-problems that are easier to solve. With such a decomposition strategy, the whole problem can be solved by optimizing the individual sub-problems independently. There are two primary decomposition directions in dealing with large-scale problems, one of which performs decomposition in the decision space by decomposing decision vectors into smaller components. In contrast, the other direction performs decomposition in the objective space by decomposing a multi-objective problem into scalar optimization sub-problems or simple single-objective sub-problems. Benders decomposition (which decomposes a problem into a master problem and a sub-problem by separating the decision variables in two groups) has been the popular decomposition method in the selected papers. Several studies adopted more unique solution algorithms, though. For example, Nested Column-and-Constraint Generation (NCCG) in \cite{fang2017optimizing}, two-stage solution algorithm in \cite{li2015bilevel}, and Multi-Criteria Decision Making (MCDM) techniques in \cite{karakoc2020social} are among the noteworthy unconventional solution methods for mathematical models.

\begin{table}[htbp]
	\centering 
	\caption{ Summary of studies developing deterministic optimization models and exact solution algorithms}
		\centering 
\begin{tabular}{ccc}
\hline
Model type     & Documents       & Solution approach                    \\
\hline
\hline
LP                         &  \cite{hu2013post}   &     Weighting method    \\
\hdashline
\multirow{2}{*}{IP}        &  \cite{fu2009optimizing} &    Direct optimization   \\
                          &  \cite{AbdelMottaleb2019exact}&  Direct optimization \\
                          \hdashline
                          
nonlinear   IP            &    \cite{he2017gas}  &    Direct optimization    \\

\hdashline
\multirow{26}{*}{MIP}      &  \cite{lee2007restoration} &     Direct optimization     \\
                          &    \cite{gong2009logic}  &  Benders decomposition     \\
                          &          \cite{matisziw2010strategic} &  Weighting method \\
                          &   \cite{cavdaroglu2013integrating}&  Direct optimization \\
                          &  \cite{ahmadi2014distribution}   &   Direct optimization    \\
                          &  \cite{sharkey2015interdependent} &  Direct optimization      \\
                          &    \cite{azad2016disruption}    &  Direct optimization  \\
                          &    \cite{wei2016quantifying}   &  Direct optimization   \\
                          &  \cite{ouyang2017mathematical} & Direct optimization   \\
                          & \cite{abbasi2017parallel}    &   Direct optimization   \\
                          & \cite{gonzalez2017efficient}  &  Direct optimization  \\
                          &  \cite{baidya2017effective} & Direct optimization \\
                          &  \cite{rong2018optimum} &   Direct optimization     \\
                          &         \cite{karakoc2019community}& Weighting method\\
                          & \cite{almoghathawi2019resilience} &   Weighting method \\
                          &   \cite{morshedlou2019restorative}&  Direct optimization  \\
                          & \cite{ge2019co} &   Direct optimization   \\
                          &   \cite{ouyang2019value}   &     Direct optimization     \\
                          & \cite{almoghathawi2019restoring}&  Weighting method \\
                          &  \cite{karakoc2020social}  &    MCDM techniques \\

                          &    \cite{chen2020failure}  &   Direct optimization    \\
                          &   \cite{sang2021resilience} &  Direct optimization      \\
                         & \cite{mohebbi2021decentralized}& Game theory and decomposition \\
                          &  \cite{almoghathawi2021exploring}&  Direct optimization  \\
                          &  \cite{rodriguez2022resilience}&Direct optimization \\

                         \hdashline

MINLP    &  \cite{li2021resilience}&    Direct optimization          \\

\hdashline
Constrainted Quadratic     &  \cite{mohagheghi2014reinforcement}&   Direct optimization          \\

\hdashline
\multirow{4}{*}{Bi-level MIP}  &     \cite{li2015bilevel} &   Two-stage solution method    \\
                          &       \cite{xiang2017coordinated}  &  Direct optimization\\
                          &       \cite{zhao2020transportation}  &  Modified active set    \\
                          &  \cite{moazeni2021formulating} &Direct optimization \\
                          \hdashline
\multirow{5}{*}{Tri-level MIP} &       \cite{fang2017optimizing} & NCCG    \\
                          &     \cite{lai2019tri}   &  Column-and-Constraint Generation\\
                          & \cite{ghorbani2020protection} & Covering decomposition\\
                          & \cite{ghorbani2021decomposition} & Integrated Benders and set-covering decomposition\\
                          & \cite{he2021tri}& NCCG\\
                          \hdashline
                          nonlinear tri-level             &    \cite{ouyang2017spatial} &   Decomposition algorithm         \\
                          \hline
\end{tabular}
\label{exact}
\end{table}

\subsubsection{Heuristic Solution Approach}
Given the intractable behavior of exact solution approaches for medium-scale and large-scale networks, heuristic methods have emerged initially to reduce the computational complexity of mathematical models and sustain the scalability of exact methods. Heuristic methods for optimization models can be implemented for the following purposes:

\begin{itemize}
    \item \textbf{Complexity reduction}: Heuristics can be adopted to provide an initial solution to help the exact methods in the exploration of a significantly reduced solution space. As the generated warm start solution by a heuristic method can be embedded into the commercial solvers, this application area is trendy in the literature. For instance,  \cite{smith2020interdependent} reduced the computational time needed to find the optimal solution for their proposed model by applying a best-response heuristic to incorporate the emergency nature of the restoration decisions.
    
    \item \textbf{Performance guarantee}: The heuristic solution algorithms can also provide an upper bound or a near-optimal solution with a performance guarantee. It should be noted that the heuristics are still considered complementary tools to exact methods in this scope. Among others, \cite{morshedlou2018work} proposed a heuristic algorithm to obtain the lower bound for the restoration crew routing problem for a relaxed formulation of the original mathematical model. The heuristic was designed to provide a feasible initial solution aligned with policies for enhancing the resilience of infrastructure networks.  
    
    \item \textbf{Approximation solution algorithms}: Despite the problem-dependent nature of heuristics, these approaches can be considered easy to implement and powerful approximate solution methods for optimization problems. In this context, heuristics are utilized as custom-tailored computational tools to find near-optimal solutions of problems directly. For example, \cite{wu2022allocation} adopted a nested heuristic algorithm to solve the challenging Defense-Attack-Recovery tri-level optimization model. The proposed heuristic could provide high-quality feasible resource allocation decisions for defensive and restorative resources to improve the strategic resilience of a power supply system.    
\end{itemize}

Table \ref{heu} summarizes the studies that adopted a heuristic solution algorithm for exact mathematical models in terms of the model type and the scope of the specific proposed heuristic method.

\begin{table}[htbp]
	\caption{ Summary of studies developing deterministic optimization models and heuristic solution algorithms}
\begin{tabular}{ccc}
\hline
Model type       & Document          & Heuristic method description                                           \\
\hline
\hline
               &  \cite{nurre2012restoring}    & Reduce computational complexity                                         \\
               &  \cite{iloglu2018integrated} & Provides an upper bound to the   original model                         \\
               & \cite{tan2019scheduling}       & Approximation algorithms with   performance guarantees                  \\
IP                & \cite{sandor2019cyber}     & Reduce the computational burden   of optimally exact solving            \\
               &  \cite{zhang2019network}    & Obtain the near-optimal solution                                        \\
                &  \cite{iloglu2020maximal}& Identifying near-optimal solutions                                    \\
               &  \cite{habib2021cascading}      &  Solve large instances of the problem \\
\hdashline
 \multirow{8}{*}{MIP}   &  \cite{kalinowski2015incremental}& Greedy heuristics with performance guarantees                         \\
               &  \cite{kasaei2016arc}& Achieve near-optimal or optimal solutions quickly                     \\
               &  \cite{akbari2017multi}& Heuristics with performance guarantees                                \\
               &    \cite{smith2020interdependent}   & Heuristics with performance guarantees                                \\
             &  \cite{morshedlou2018work} & Provide an upper bound on the solution                                \\
               & \cite{morshedlou2021heuristic} & Generates reliable solution for large-scale problems                  \\
              &   \cite{kong2021optimizing}   & Identifying near-optimal   solutions                                    \\
              &       \cite{li2018collaborative}   & A greedy heuristic to solve the optimization problem                  \\
\hdashline
 \multirow{2}{*}{Tri-level MIP} &   \cite{fakhry2022tri}   & A trade-off between solution quality and computational time           \\
 &   \cite{wu2022allocation}       & Variable neighborhood heuristic to solve the model \\       \hline           
\end{tabular}
\label{heu}
\end{table}

\subsubsection{Evolutionary Solution Approach}
As the size of models increases (in terms of decision variables and constraints) and the exact and heuristics solution algorithms face difficulties, adopting evolutionary algorithms to find near-optimal solutions seems to be an emerging direction in the literature. The evolutionary algorithms in the selected papers can be divided based on the nature of the objective function in the model. For single objective models, for example, \cite{maya2011grasp} proposed an evolutionary solution algorithm based on the GRASP and Variable Neighborhood Search (VNS) metaheuristics for the resource allocation problem in the recovery of a rural road network after a natural or man-made disaster. The solution method was applicable to a range of problems, from small and medium-size instances to a large real-life motivated instances. As for multi-objective models, \cite{fang2014optimal} developed an Non-dominated Sorting Genetic Algorithm (NSGA-II), the multi-objective extent ion of the GA algorithm, to tackle the power transmission network resilience problem in the presence of cascading failures with limited investment costs. The problem was formulated as a multi-objective optimization problem to find a trade-off between cost and the system's vulnerability to cascading failures. Table \ref{ev} presents a summary of evolutionary algorithms for deterministic models, including Dynamic Programming (DP), in the context of CPS regarding the type of model and the adopted evolutionary solution algorithm. The summarized information for evolutionary methods for deterministic models shows that using GA and NSGA-II are the most popular solution methods due to their flexibility in dealing with various models and technical assumptions. However, a change of direction towards more powerful algorithms such as GWO in \cite{kumar2022hybrid} and ACO (which is shown to be effective in dealing with network-based problems) in \cite{li2019joint} is noticeable in this area.    

\begin{table}[htbp]
\small
	\caption{ Summary of studies developing deterministic optimization models and evolutionary solution algorithms}
\begin{tabular}{ccc}
\hline
Model   type          & Document                & Evolutionary solution approach                               \\
\hline
\hline
Constrained  MIP &  \cite{kumar2022hybrid}       & Grey Wolf Optimization (GWO)                                 \\
\hdashline
\multirow{2}{*}{DP}             &   \cite{zhang2018restoration}         & GA                                                           \\
            & \cite{duque2016network}           & Iterated Greedy-Randomized Constructive Procedure  \\
\hdashline
\multirow{10}{*}{IP}               &     \cite{orabi2009optimizing}        & Multi-objective genetic   algorithm                          \\
                                   &  \cite{orabi2010optimizing}          & NSGA-II                                                      \\
                                   &     \cite{sen2011optimized}          & GA                                                           \\
                                   &     \cite{maya2011grasp}       & GRASP and VNS                                                \\
                                   &    \cite{vodak2018modified}        & Ant Colony Optimization (ACO)                                     \\
                                   &    \cite{zhang2018optimal}        & Particle Swarm Optimization (PSO)                                  \\
                                   &     \cite{li2019joint}          & ACO                                     \\
                                   &     \cite{sharma2020regional}       & GA                                                           \\
                                   &    \cite{ghiasi2021resiliency}         & GWO                                  \\
                                   &    \cite{zhu2021integration}          & Probabilistic Solution Discovery Algorithm                               \\
                                   \hdashline
\multirow{4}{*}{MIP}               &  \cite{pramudita2014model} & TS                                                  \\
                                   &  \cite{ouyang2015resilience}        & GA                                                           \\
                                   &     \cite{zuloaga2020interdependent}       & GA                                                           \\
                                   &    \cite{poudel2020generalized}         & GA                                                           \\
                                   \hdashline
\multirow{2}{*}{Nonlinear   LP}    &       \cite{fang2014optimal}       & NSGA-II                                                      \\
                                   &    \cite{fang2015optimization}          & Non-dominated Sorting Binary Differential Evolution \\
                                   \hdashline
Nonlinear  MIP                    &      \cite{song2016intelligent}        & Orthogonal GA                                                \\
\hdashline
Nonlinear bi-level MIP        &  \cite{wang2013combined}        & SA  \\                               \hline                         
\end{tabular}
\label{ev}
\end{table}
\newpage

\subsection{Stochastic Formulation}
\begin{figure}[htbp]
    \centering
\begin{center}
\begin{tikzpicture}[scale=1.3]
\begin{axis}[
    ybar stacked,
    axis lines=left,
    enlarge x limits=0.03,
    ymin=0,
    x tick label style={/pgf/number format/1000 sep=,rotate=90,anchor=east, font=\footnotesize},
    y tick label style={/pgf/number format/1000 sep=,anchor=east, font=\footnotesize},
    legend style={draw=none,at={(0.22,0.97)},anchor=north, font=\footnotesize},
    height=6 cm,
	width=12 cm,
    x label style={draw=none,at={(axis description cs:0.5,-0.10)},anchor=north, font=\footnotesize},
	xlabel={Year},
	y label style={at={(0.05,0.5)}},
	ylabel style={font=\footnotesize},
    ylabel={No. of documents},
    symbolic x coords={2007, 2008, 2009, 2010, 2011,2012,2013,2014,2015,2016,2017,2018,2019,2020,2021,2022
		},
    xtick=data,
     ymin = 0, ymax = 7,
    ]
\addplot+[style={RawSienna},ybar] plot coordinates {

(2007,0)
(2008,0)
(2009,0)
(2010,0)
(2011,0)
(2012,0)
(2013,1)
(2014,0)
(2015,2)
(2016,1)
(2017,1)
(2018,4)
(2019,4)
(2020,3)
(2021,1)
(2022,1)
};
\addplot+[style={darkgray},ybar] plot coordinates {
(2007,0)
(2008,0)
(2009,1)
(2010,0)
(2011,0)
(2012,0)
(2013,0)
(2014,0)
(2015,1)
(2016,0)
(2017,0)
(2018,1)
(2019,1)
(2020,1)
(2021,2)
(2022,0)
};
\addplot+[style={},ybar] plot coordinates {
(2007,1)
(2008,1)
(2009,0)
(2010,0)
(2011,0)
(2012,0)
(2013,0)
(2014,0)
(2015,1)
(2016,1)
(2017,0)
(2018,0)
(2019,1)
(2020,0)
(2021,1)
(2022,0)
};
\legend{\strut Exact solution algorithm, \strut Heuristic solution algorithm, \strut Evolutionary solution algorithm}
\end{axis}
\end{tikzpicture}
\caption{Number of papers per year and per each type of solution method for stochastic mathematical model.}
\label{sto}
\end{center}
\end{figure}
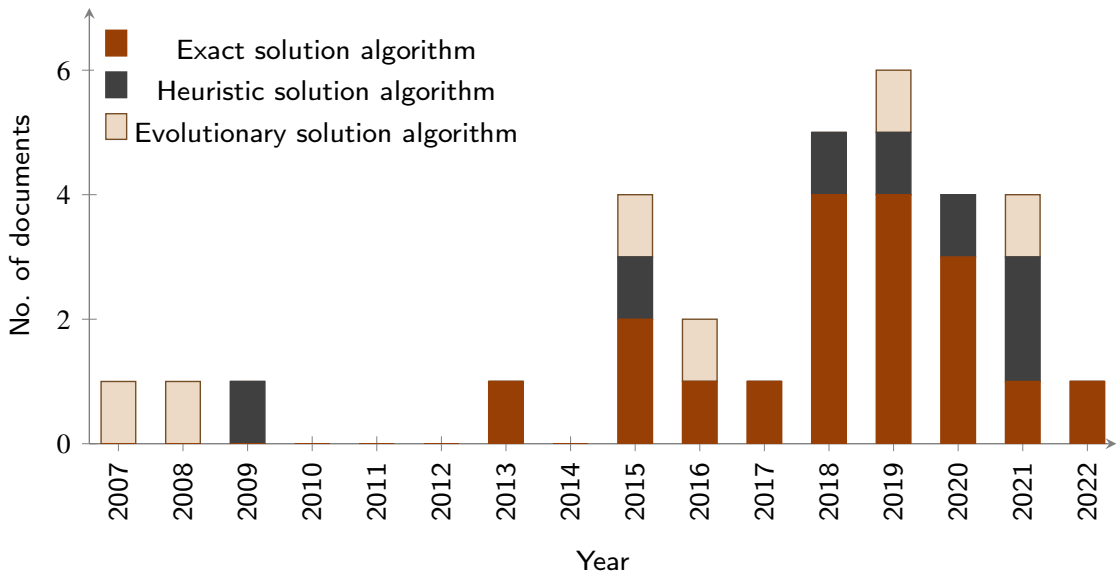

In this section, we focus on the studies that developed optimization solution methods for stochastic mathematical models and classify those works based on their solution approach. We present a comprehensive overview of the solution algorithms in each sub-section and some representative studies. We also explore how selected papers incorporated the concept of uncertainty in their modeling step. Fig. \ref{sto} shows the frequency of each solution approach for stochastic formulation based on the publication date of documents.

\subsubsection{Exact Solution Approach}
Given the high level of uncertainty in the post-disaster environment (e.g., the aftershock in the case of earthquakes or secondary disruptions), increasing attention has been dedicated to developing stochastic optimization models for cyber-physical networks. Focusing on the solution algorithm, decomposition methods such as Benders decomposition and accelerated L-shape algorithms are popular options for two-stage stochastic models. For instance, \cite{arab2015stochastic} formulated a proactive resource allocation model to repair and restore potential damages to the power system infrastructure as a two-stage stochastic integer model with complete recourse. They adopted Benders decomposition to handle the computational hurdle of their model's large-scale mixed-integer equivalence counterpart. Several studies such as \cite{kavousi2018stochastic} utilized simulation-based methods as a complementary tool for the optimization algorithms in an uncertain environment. \cite{zhao2017hybrid} is one unique work in this class by developing a Hybrid Hidden Markov Model (HMM) to optimize dispatch strategy in post-disruption stage to maximize system resilience and ensure the acceptable system functionality level. Table \ref{stex} summarizes the studies with an exact solution algorithm for the stochastic model and associated sources of uncertainty in the stochastic formulation. The information in this table shows that repair time, demand, and failure are the most common source of randomness in the developed models.

\begin{table}[htbp]
\footnotesize
\caption{ Summary of stochastic programming studies}
\begin{tabular}{cccc}
\hline
Model type                  & Document               & Sources of uncertainty                                  & Solution approach        
\\
\hline
HMM       & \cite{zhao2017hybrid}       & Capacity of supply                                                       & Optimized   dispatch strategy                          \\
\hdashline

LP                          & \cite{holden2013network} & Functionality & Monte Carlo simulation                                 \\
\hdashline
MILP                        & \cite{wu2021risk} & Repair time and resources                         & Decomposition                                          \\
\hdashline
MIP                         & \cite{gonzalez2016interdependent}    & Demands, cost, and functionality      & Simulation+Optimization                                \\
\hdashline
Nonlinear MIP              & \cite{zhang2018resilience}       & Disruption                                                        & Probabilistic solution discovery   \\
\hdashline
\multirow{12}{*}{Two-stage  MIP} & \cite{arab2015stochastic}    & Damage   state                            & Benders decomposition                                  \\
                            &  \cite{arif2018optimizing}     & Repair time \& demand                                   &  Progressive Hedging                    \\
                            & \cite{kavousi2018stochastic} & Travel   time                                           & Simulation                                             \\
                            &\cite{beheshtian2018climate}   & Failure                                                 & Simulation                                             \\
                            & \cite{fang2019optimum}& Repair   crews and Restoration time                      & Tailored Benders decomposition                          \\
                            &\cite{fang2019adaptive}     & Failure                                                 & Nested decomposition                                   \\
                            &  \cite{sanci2019integrating}  & Demand, supply, and network availability & Sample average approximation \\
                            & \cite{gomez2019optimization}    & Disruption                                              & Decomposition techniques                               \\
                            & \cite{bhuiyan2020stochastic}     & Survival probability of a component                     & Accelerated L-shaped                                   \\
                            &  \cite{wang2020resilience}       & Disruption                                              & L-shaped method                                        \\
                            &\cite{abessi2020new}     & Load consumption and demand response                   & Simulation                                             \\
                            & \cite{alkhaleel2022risk} & Repair time \& travel time& Benders decomposition  \\
                            \hline
\end{tabular}
\label{stex}
\end{table}

\subsubsection{Heuristic Solution Approach}
A few papers adopted heuristic methods to solve the computationally challenging stochastic models. Similar to deterministic models, heuristics have been developed for two various purposes. First, the heuristic methods aimed to reduce the computational complexity. For instance, \cite{ccelik2015post, yang2018towards} used a heuristic scheme to simplify the core optimization problem, Stochastic Dynamic Programming (SDP), of link recovery selection in a transportation network to improve the resilience. Second, solving the model approximately was sought in \cite{ulusan2021approximate} by developing an Approximate Dynamic Programming (ADP) approach to solve the stochastic road network recovery problem heuristically. Table \ref{hest} summarizes the studies with a heuristic solution algorithm for the stochastic model and associated sources of uncertainty in the stochastic formulation. In this category, repair time, and failure (or equivalently disruption) are the most common source of uncertainty.

\begin{table}[htbp]
\caption{ Summary of studies developing stochastic optimization models and heuristic solution algorithms}
\begin{tabular}{cccc}
\hline
Model type                       & Document         & Sources of uncertainty   & Heuristic scope                                                                          \\
\hline
\hline
\multirow{2}{*}{SDP} &\cite{ccelik2015post}     & Debris   amount     & Reduce the computational time                            \\
&   \cite{yang2018towards}   & Recovery time       & Reduce computational complexity                                                          \\
\hdashline
\multirow{2}{*}{MDP}  &  \cite{nozhati2020optimal}   & Repair time         & Rollout technique to find   optimal policy for an MDP                                    \\
 &  \cite{ulusan2021approximate}     &Demand       &  Solve the model approximately       \\
 \hdashline
\multirow{3}{*}{MIP}                           &   \cite{liu2021heuristic}    & Disruption severity & Solve the problem efficiently \\
                            &  \cite{che2009lagrangian}    & Disruption duration & Find near optimal solutions           \\
                            & \cite{fang2019resilient}    & Disruption          & Variable Neighborhood Search For   solving the problem\\
\hline
\end{tabular}
\label{hest}
\end{table}

\subsubsection{Evolutionary Solution Approach}
The literature for this solution approach is minimal, and only a handful of papers developed evolutionary algorithms for stochastic problems for CPS affected by extreme disruptive events. For example, \cite{liu2021hierarchical} employed NSGA-II to solve the multi-objective stochastic optimization problem for simultaneously minimizing cost and maximizing the resilience of an infrastructure network. Table \ref{evst} summarizes the studies with an evolutionary solution algorithm for the stochastic model and associated sources of uncertainty in the stochastic formulation. Similar to other sub-sections of stochastic programming, repair time is the most common source of uncertainty.

\begin{table}[htbp]
\caption{ Summary of studies developing stochastic optimization models and heuristic solution algorithms}
\begin{tabular}{cccc}
\hline
Model type                       & Document       & Source of uncertainty                           & Solution method                   \\
\hline
\hline
\multirow{2}{*}{IP}                              &\cite{xu2007optimizing}     & Restoration time                                & GA                                \\
                               &  \cite{zou2019enhancing} & Disruptions,   traffic demands, and repair cost & Binary PSO      
                                \\
                                \hdashline
\multirow{2}{*}{MIP}                             & \cite{furuta2008optimal} & Restoration time                                & GA considering uncertainty (GACU)  \\
                              &  \cite{zhang2016two}  & Restoration time                                & NSGA-II                           \\
                               \hdashline
SDP&   \cite{liu2021hierarchical} & Restoration time and horizon                     & NSGA-II \\
\hline
\end{tabular}
\label{evst}
\end{table}

\subsection{Machine Learning-Enabled Solution Approaches}
As CPS networks become more comprehensive, the optimization problems associated with the city-scale systems exhibit computational difficulties. From another point of view, ML has been regarded as a promising research direction to address the challenges mentioned above in conventional optimization algorithms \citep{mirshekarian2018machine}. ML can improve traditional solution algorithms for a range of optimization problems from two aspects. First, assuming that there is a piece of expert knowledge (\textit{a priori} knowledge) about the optimization algorithm, learning can be used to build approximations to alleviate the computational effort \citep{lodi2020learning}. Second, in the absence of expert knowledge, the goal of a learning procedure is to explore the decision (or objective) space of these problems to lead the search process to the optimal (or near-optimal) performing behavior. Consequently, these intelligent decisions improve the performance of optimization algorithms in terms of solution quality, convergence rate, and robustness \citep{karimi2020learning}. 

There are a few studies in recent years leveraging ML algorithms and data-driven methods to augment the optimization algorithms for CPS. For example, \cite{alemzadeh2020resource} developed a data-driven approach based on Artificial Neural Networks (ANNs) to estimate the optimal restoration sequence for disrupted infrastructure networks after natural disasters such as earthquakes. \cite{zhao2020data} also proposed a data-driven framework for solving the optimal two-stage stochastic optimization problem to enhance the reliability of planning and operation of a power-gas network against seismic attacks. \cite{yang2021real} combined the pre-training of an artificial intelligent agent and optimization techniques to accelerate calculating the optimum recovery path in the case of electrical faults in an interdependent network. Most recently, \cite{aslani2022learn} embedded a learning component to learn from the obtained non-dominated solutions of a decomposition-based evolutionary algorithm periodically and consequently improve the performance of the optimization algorithms for multi-objective models in restoration of large-scale infrastructure networks.

\section{Literature Analysis} \label{analyze}
This section provides a comprehensive set of analyses from bibliography, modeling, and scalability points of view. In each sub-section, we highlight the key observations and patterns to shed light on the current trends and significant gaps in the optimization solution algorithms for CPS.   
\subsection{Bibliographic Analysis}
As the first step, we analyzed the pool of literature review to obtain key insights from the current patterns. Specifically, we focus on two bibliographic analyses, keyword clustering and bibliographic coupling. The analyses were carried out in VOSviewer™ software, which is developed based on  \cite{eck2007vos} to present a new method for visualizing similarities between objects. As a preliminary step, we explored all the selected papers in \textbf{Web of Science (WOS)} core collection. Among the chosen articles, 3 studies were not indexed in this database, which is less than 3\% of the pool for review. In the next step, we imported all the bibliographic information of papers into the VOSviewer™. We then conducted two clustering procedures as follows: 

\begin{itemize}
    \item \textbf{Keyword clustering}: In this part, we focused on the keywords of selected papers. The output of this analysis is visualized in Fig. \ref{key}. Following a series of trial and error, we set the parameter $n=7$ as the threshold for the overall occurrence of a keyword to be included in the clustering process. In this figure, the label size and the circle indicate the frequency of the associated keyword. Therefore, the keywords with higher frequency are depicted in a larger font and are more legible. In addition, the color of a keyword mirrors the assigned cluster to the keyword. Finally, the distances between keywords show the relatedness in the literature regarding co-occurrence links. In other words, closer terms indicate closer relatedness and more dense lines represent the strongest co-occurrences. The output shows that there is a balance among methodology-focused keywords (e.g., \textit{optimization, model, and design} clustered in green) and application-based keywords (e.g., \textit{risk, resilience, and reliability} clustered in red). In addition, the term \textit{power} is a frequent keyword, showing that power networks are among the most investigated CPS. This observation is consistent with the further analysis conducted in section \ref{application}.

\begin{figure}[htbp]
	\centering
	\includegraphics[scale=0.22]{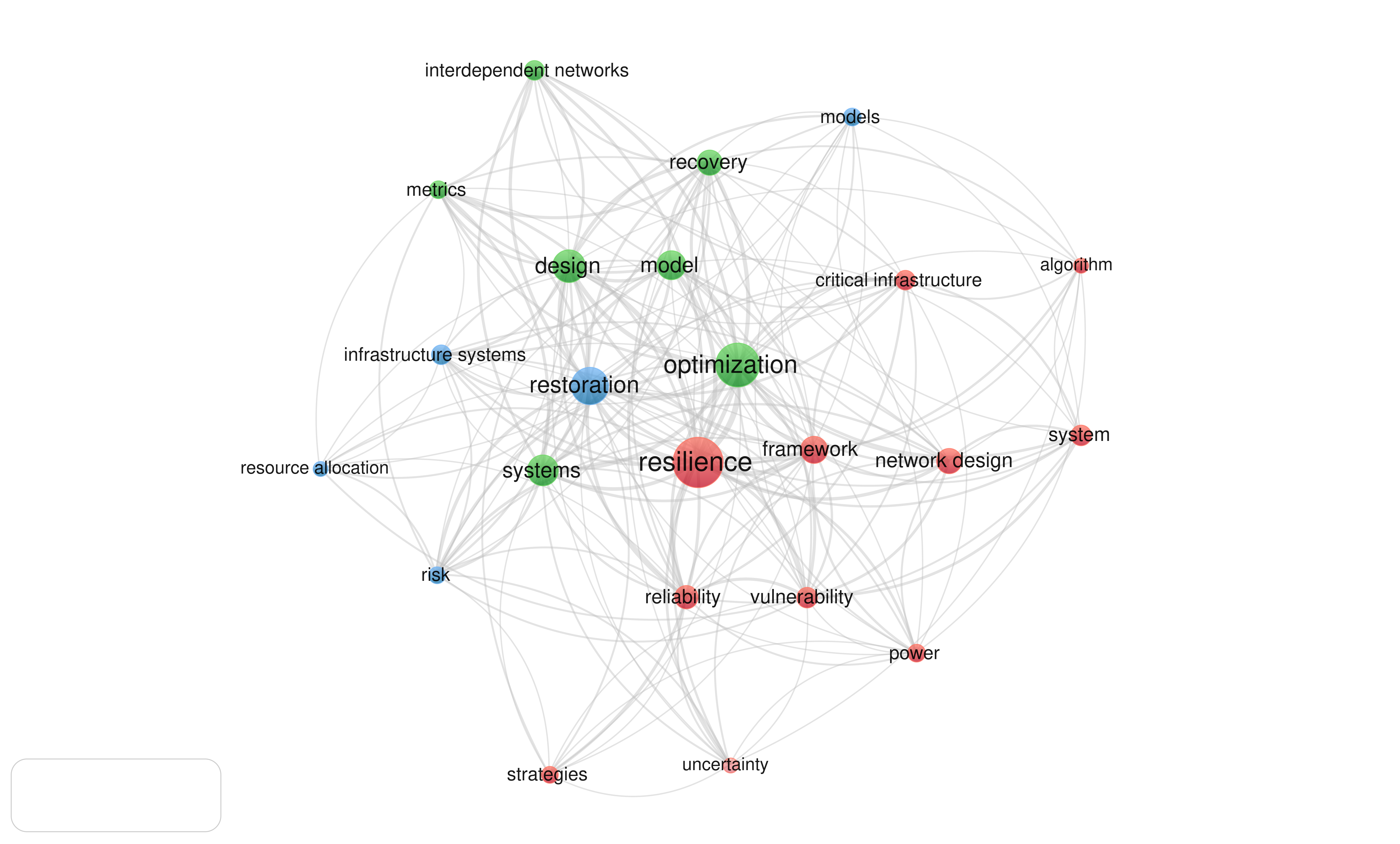}
	\centering
	\caption{Keywords clustering based on co-occurrence visualized by VOSviewer™ (n = 7)}
	\label{key}
\end{figure}

\item \textbf{Bibliographic coupling}: In this analysis, we focused on the peer-reviewed journals (imported from WOS) published the selected set of papers in the review. Bibliographic coupling is defined when two works cite a third common work in their bibliographies. This index is an indication that a probability exists that the two coupled research papers approached a related subject matter. The \textit{coupling strength} of two papers is higher when they have more citations to other papers they share (reflected with a thicker line in the final figure). The output of this analysis is visualized in Fig. \ref{journal}. Similar to the other experiment, we conducted a parameter tuning method and set $n=2$ as the threshold for inclusion, meaning there is a coupling if two papers at least cited two common references. Each color in this figure represents one cluster of journals, and the connections show that there is a strong mutual coupling among the works published in Operations Research (OR) methodological journals (denoted by green nodes) such as \textit{European Journal of Operational Research (EJOR)} and \textit{Computers and Industrial Engineering (CAIE)} and application-focused journals (blue and red nodes) such as \textit{Sustainable cities and society} and \textit{safety science}. In addition, we can observe a strong clustering among the journals related to \textit{transportation} (green cluster) and \textit{power} (red cluster) infrastructure systems.   
\begin{figure}[htbp]
	\centering
	\includegraphics[scale=0.17]{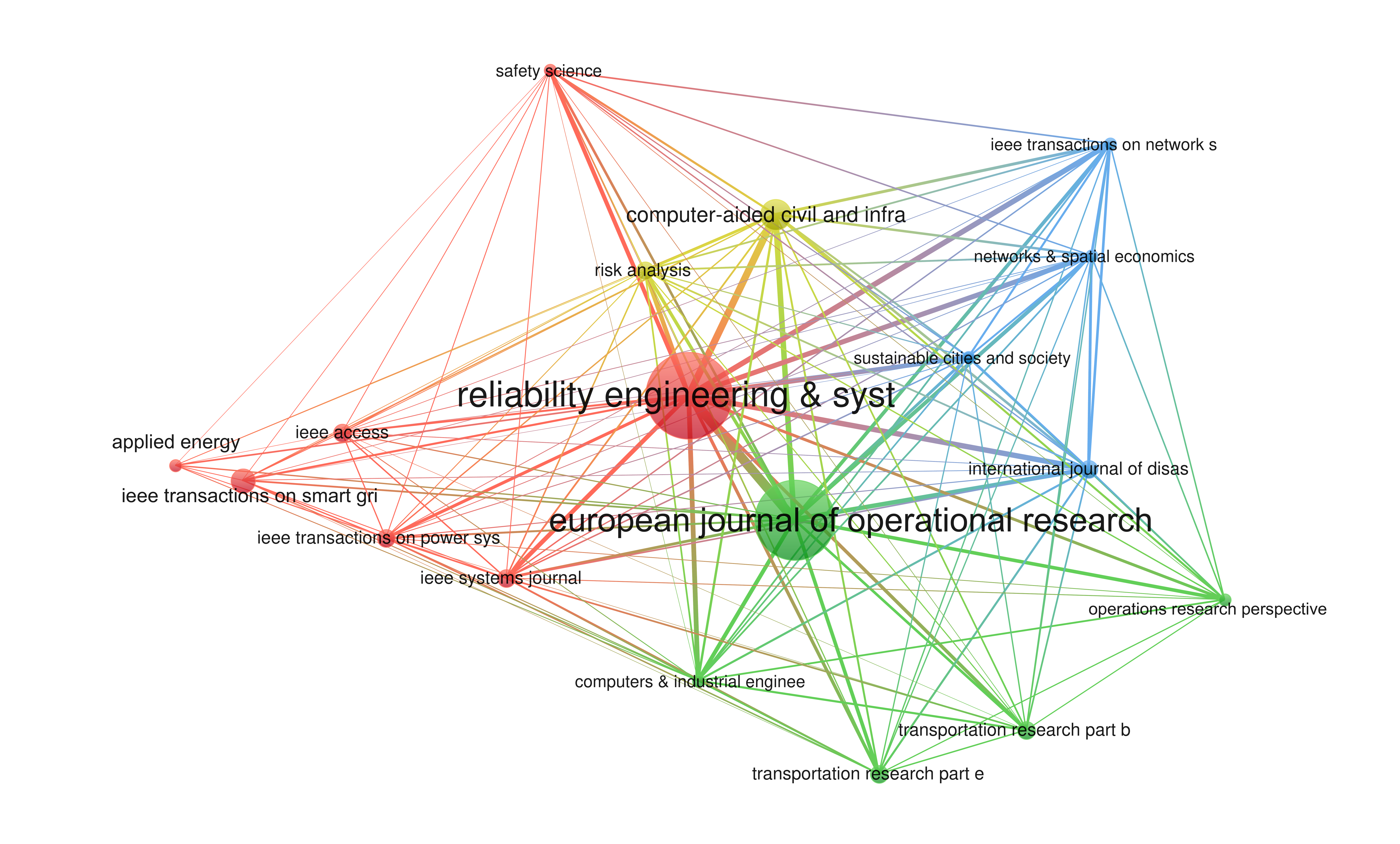}
	\centering
	\caption{Bibliographic coupling visualized by VOSviewer™ (n = 2)}
	\label{journal}
\end{figure}
\end{itemize}

\subsection{Decision Variables}
Decision variables are the major component of mathematical modeling as they form the model and drive the set of constraints. In a general classification, decision variables can be continuous (accepting values from zero to an upper bound) or integer (taking only integer values). Binary decision variables are a specific class of integer category by taking only 0 and 1 values. In the literature of CPS, both types of decision variables are present in high density. Table \ref{dv} shows a summary of decision variables with one representative study for each variable type.

Continuous decision variables have been utilized to capture the specific features of infrastructure networks. For instance, the flow of an arc (defined as a physical connection among two nodes like a water pipe) is a standard decision variable and is present in several works such as \cite{lee2007restoration} and \cite{sharkey2015interdependent}. A more flow-related specific decision variable is the cumulative flow of demand nodes (supply nodes should provide the commodity to meet demand requirements), which is employed in \cite{morshedlou2019restorative}. The excess supply in \cite{gonzalez2017efficient} and unmet demand in \cite{almoghathawi2019restoring} are the other popular continuous decision variables in this application area.  

On the other hand, binary decision variables are more diverse. As many decision in the restoration phase is of a binary nature, this type of variables can adequately capture the complexity of decision-making procedure. For instance, the arc traversal decision variable in \cite{kasaei2016arc}, protection in \cite{zhu2021integration}, attacker decision in \cite{ghorbani2021decomposition}, and resource allocation in \cite{nurre2012restoring} reflect a similar decision about the incorporation of a component in the restoration from different perspectives. Likewise, the network design decision variables in \cite{iloglu2018integrated} and temporary arc installment in \cite{cavdaroglu2013integrating} have a similar nature in the physical extension of a network permanently and temporarily, respectively. Finally, the scheduling decision variables in a variety of works such as \cite{karakoc2019community}, child-parent interdependency in \cite{sharkey2015interdependent}, and functionality in \cite{matisziw2010strategic} aim to capture more sophisticated aspects of recovery planning.

\begin{table}[]
\caption{ Summary of decision variables in optimization models}
\begin{tabular}{c|ll}
\hline
Type &Decision variable                & Selected document      \\
\hline
{\multirow{4}{*}{\rotatebox[origin=c]{90}{Continuous}}}&Arc flow                     & \cite{lee2007restoration}                   \\
&Cumulative   flow to demand node &\cite{morshedlou2019restorative}  \\
&Excess supply                   & \cite{gonzalez2017efficient}     \\
&Unmet demand                    & \cite{almoghathawi2019restoring}\\
\hdashline
{\multirow{9}{*}{\rotatebox[origin=c]{90}{Binary}}}&
Arc traversal                     & \cite{kasaei2016arc}     \\
&Attacker decision                 &\cite{ghorbani2021decomposition}  \\
&Child-parent interdependency      &\cite{sharkey2015interdependent}             \\
&Functionality                    & \cite{matisziw2010strategic}               \\
&Network design                    &\cite{iloglu2018integrated}    \\
&Protection decision               & \cite{zhu2021integration}\\
&Resource allocation              & \cite{nurre2012restoring}                  \\
&Scheduling                        & \cite{karakoc2019community}         \\
&Temporary arc   installment        & \cite{cavdaroglu2013integrating}             \\

\hline
\end{tabular}
\label{dv}
\end{table}

\subsection{Objective Functions} \label{object}
\begin{figure}[htbp]
	\centering
	\includegraphics[scale=0.55]{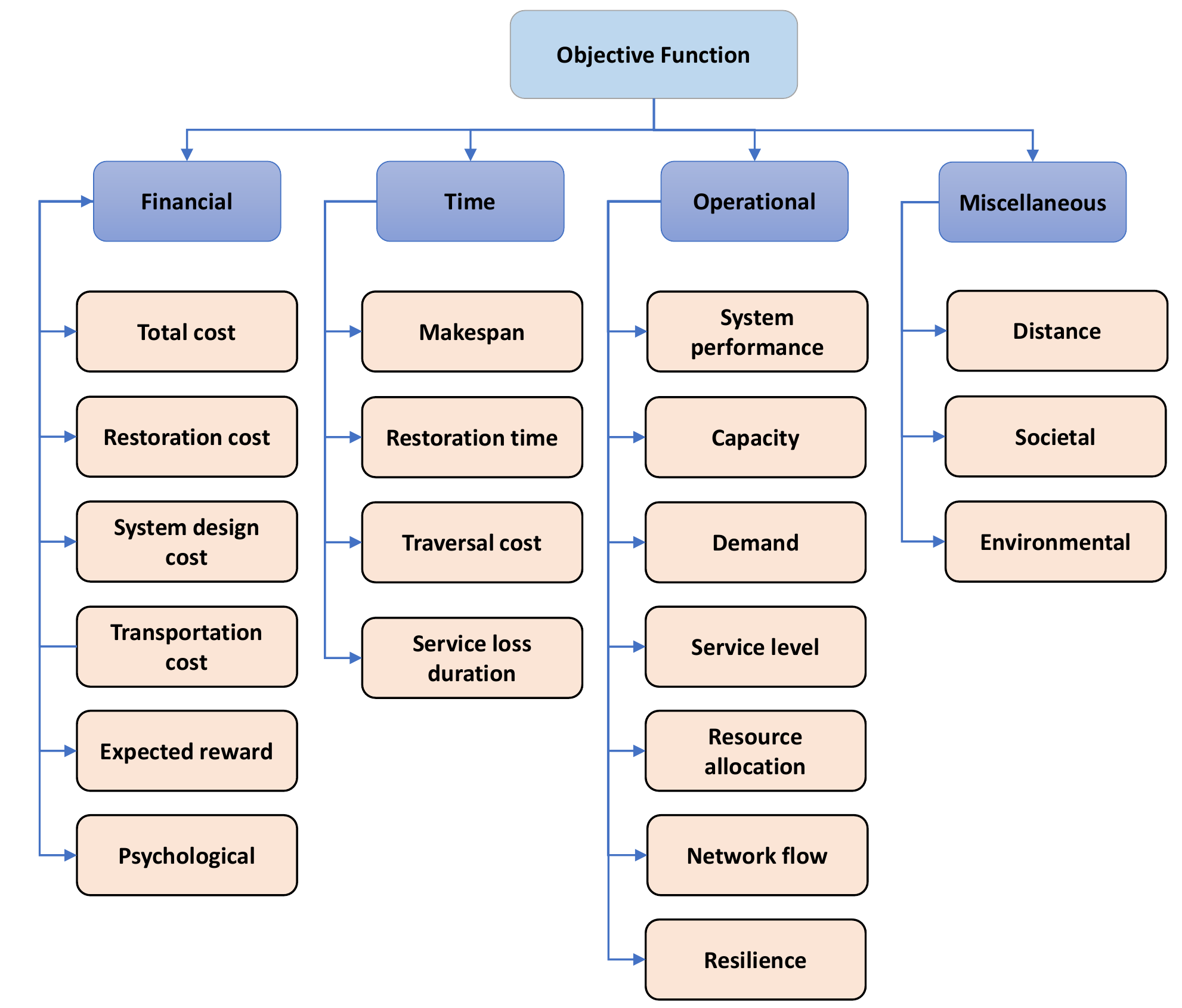}
	\centering
	\caption{Objective function classification}
	\label{mc}
\end{figure}

The studies in the literature defined various objective functions to capture significant aspects of a disrupted CPS. This section briefly provides a classification of objective functions and their characteristics. Table \ref{obj} presents a classification of objective functions for all papers reviewed in this study. Fig. \ref{mc} also shows the classification of objective functions present in the literature.  

The financial part of restoration plays a crucial role in both low-level (i.e., the component level) and high-level (i.e., the network level) decision-making layers. On the high level, the network administrators aim to minimize the total cost of the CPS system. The term \textit{total cost} in the literature referred to a summation of different cost components forming a monolithic cost figure. For example, \cite{almoghathawi2019resilience} defined the system cost as the objective function, including restoration cost, flow cost, and disruption cost (i.e., unmet demand). The system design cost is another example of a high-level cost function considered in the work of \cite{zhang2018resilience}. On the other hand, several papers focused on the low-level cost, such as transportation cost in \cite{bhuiyan2020stochastic}, restoration (or recovery) cost in \cite{smith2020interdependent}. Finally, there are some unique financial objective functions, such as maximizing the expected reward in \cite{ulusan2021approximate} and minimizing the psychological costs during the restoration phase in \cite{hu2013post}.  

As the restoration problem is essentially a scheduling problem, time-related objective functions are also a common option in the modeling phase. Among others, \cite{zhang2016two} defined the makespan (the time difference between the start of restoration and finishing the last assigned repair) to expedite the recovery of an affected system to the standard performance level. Minimizing the restoration time is also investigated in \cite{vodak2018modified}, and \cite{orabi2010optimizing} attempted to reduce the duration of service loss. Finally, several papers such as \cite{kasaei2016arc,akbari2017multi, zhang2019network}, which focused on transportation infrastructures, defined specific objective functions to minimize the traversal cost of disrupted arcs or the longest walk (transformed into an equivalent time component). 

Operational objective functions are the next category, in which the focus of the optimization is devoted to the system-level performance. For instance, \cite{fang2017optimizing} formulated the defender-attacker-defender optimization model for the resilience against intentional attacks to maximize the \textit{system performance} quantified as the normalized total satisfied demand level. Maximizing \textit{operational capacity} of infrastructure systems is another relevant objective function defined in \cite{ahmadi2014distribution} for an interdependent network, including a hospital, to save more lives in the case of a natural disaster. Besides, the total met demand can also reflect the performance level of an infrastructure network, which is captured in several works such as maximizing the total met demand \cite{liu2021heuristic} and minimizing the unmet demand in \cite{ghorbani2020protection}. A share of literature body also defined objective functions solely based on the network flow ( e.g., \cite{matisziw2010strategic}).

The term \textit{resilience} is also defined as the objective function of many papers from different perspectives. While this performance metric is essentially an operational measure, this objective function has been approached from distinct points of view in the CPS literature. Resilience, as a dynamic property of complex systems, is a performance measure for characterizing the behavior of CPS networks in the presence of disruptive events. Resilience generally implies the ability of a system to return to normal condition after a disturbance \citep{ouyang2015resilience}. However, this term has been defined differently across various disciplines, including ecology, emergency management, engineering, and national security. In other words, there is no consensus on the definition of this term. For example, the \textit{National Research Council} defines resilience as "The ability to prepare and plan for, absorb, recover from, or more successfully adapt to actual or potential adverse events", which is more focused on the emergency management discipline. As another definition, \textit{Social–ecological} viewpoint further emphasizes transformability, learning, and innovation as crucial aspects of complex systems resilience. The \textit{community resilience} is also defined as "the ability of a community to prepare for anticipated hazards, adapt to changing conditions, and withstand and recover rapidly from disruptions" (see \cite{sharkey2021search} for a comprehensive review). 

There are a few studies that defined distance-based objective functions such as cumulative multiple coverage of emergency demand
in \cite{iloglu2018integrated} and the longest walk in \citep{akbari2017multi}. Societal goals are also present in terms of Human well-being level in \cite{yang2021real} and the social impact of failed component restoration in \cite{aslani2022learn}. Finally, \cite{aslani2022learn} is the only study that embedded an environmental-related objective function (in terms of $CO_2$ footprint of physical components) in a multi-objective model. We classified these rare objective functions under miscellaneous category in our analysis visualized in Fig. \ref{mc}.

\begin{sidewaystable}[p]
\caption{ Classification of objective function in optimization models}
\footnotesize
\begin{tabular}{c|c|l}
\hline
Main category                & Sub-category           & Document         \\ \hline \hline
\multirow{11}{*}{Financial}        & \multirow{6}{*}{Total cost }              &    \cite{lee2007restoration}    \cite{fu2009optimizing} \cite{gong2009logic} \cite{matisziw2010strategic} \cite{cavdaroglu2013integrating}   \cite{azad2016disruption} \cite{wei2016quantifying}   \\

     &           &     \cite{gonzalez2017efficient} \cite{he2017gas} \cite{karakoc2019community}  \cite{almoghathawi2019resilience} \cite{karakoc2020social} \cite{chen2020failure} \\
     
          &             & \cite{almoghathawi2021exploring} \cite{rodriguez2022resilience} \cite{smith2020interdependent} \cite{li2018collaborative} \cite{pramudita2014model} \cite{fang2014optimal} \\
          
                    &             &  \cite{fang2015optimization} \cite{sharma2020regional} \cite{ghiasi2021resiliency} \cite{zuloaga2020interdependent} \cite{song2016intelligent} \cite{holden2013network} \cite{gonzalez2016interdependent} \\
                    
                             &             &   \cite{arab2015stochastic} \cite{sanci2019integrating} \cite{gomez2019optimization}  \cite{che2009lagrangian} \cite{fang2019resilient} \cite{liu2021hierarchical} \cite{zhang2016two} \\
                             
                                                          &              &    \cite{alemzadeh2020resource} \cite{zhao2020data}\\
\cline{2-3}
                             & Financial loss         &   \cite{AbdelMottaleb2019exact}               \\
                             & Recovery cost          &    \cite{orabi2010optimizing}     \cite{wang2013combined} \cite{aslani2022learn}         \\
                             & System design cost          &      \cite{zhang2018resilience}            \\
                             & Transportation cost       &     \cite{bhuiyan2020stochastic}             \\
                             & Expected Reward                 &       \cite{nozhati2020optimal} \cite{ulusan2021approximate}           \\
                             & Psychological cost         &   \cite{hu2013post}               \\
                               \hdashline
\multirow{4}{*}{Time  }                       & Makespan               &   \cite{gong2009logic} \cite{zhang2016two}               \\
                             & Restoration time           &  \cite{tan2019scheduling} \cite{vodak2018modified}               \\
                             & Traversal time              &      \cite{kasaei2016arc} \cite{zhang2019network}  \cite{maya2011grasp}          \\
                             & Service loss  duration         &  \cite{xu2007optimizing}                 \\
                               \hdashline
\multirow{14}{*}{Operational} & {System performance}&  \cite{fang2017optimizing} \cite{orabi2009optimizing} \cite{duque2016network} \cite{jena2021design} \cite{beheshtian2018climate}   \cite{furuta2008optimal} \cite{zou2019enhancing} \\
                             & Capacity               &    \cite{ahmadi2014distribution}  \cite{AbdelMottaleb2019exact}            \\
                             \cline{2-3}
                             & \multirow{2}{*}{Demand}                 &   \cite{he2017gas}  \cite{ghorbani2020protection}     \cite{li2021resilience} \cite{iloglu2020maximal}    \cite{sen2011optimized} \cite{li2019joint} \cite{liu2021heuristic}  \\
                                 &          &  \cite{mohebbi2021decentralized} \cite{aslani2022learn}\\
                                 \cline{2-3}
                             & \multirow{2}{*}{Service Level}         &                 \cite{mohagheghi2014reinforcement} \cite{sharkey2015interdependent} \cite{xiang2017coordinated} \cite{baidya2017effective} \cite{ge2019co} \cite{lai2019tri} \cite{he2021tri} \cite{sandor2019cyber}  \\
                                &          &                  \cite{habib2021cascading} \cite{poudel2020generalized} \cite{arif2018optimizing} \cite{abessi2020new} \cite{loggins2015rapid} \\
                                \cline{2-3}
                             & {Resource   allocation}& \cite{li2015bilevel}\\
                                                        &Network Flow   &        \cite{matisziw2010strategic}     \cite{rong2018optimum} \cite{nurre2012restoring} \cite{kalinowski2015incremental}  \cite{ccelik2015post}     \\   
                             \cline{2-3}
                                                 &  \multirow{6}{*}{Resilience}   &     \cite{ouyang2017mathematical} \cite{abbasi2017parallel}   \cite{ouyang2017spatial}  \cite{karakoc2019community}   \cite{almoghathawi2019resilience}   \cite{morshedlou2019restorative}  \\ 

                 &                        &  \cite{ouyang2019value} \cite{zhao2020transportation} \cite{karakoc2020social} \cite{almoghathawi2019restoring} \cite{moazeni2021formulating}  \\ 

                  &                        &   \cite{ghorbani2021decomposition} \cite{sang2021resilience} \cite{almoghathawi2021exploring} \cite{morshedlou2018work} \cite{morshedlou2021heuristic} \cite{kong2021optimizing}  \\ 

                 &                        &   \cite{huang2022optimization} \cite{fakhry2022tri} \cite{wu2022allocation} \cite{fang2014optimal} \cite{ouyang2015resilience} \cite{fang2015optimization} \cite{zhang2018optimal}  \\

                  &                        &  \cite{sharma2020regional} \cite{zhang2018restoration} \cite{zhu2021integration} \cite{zhao2017hybrid} \cite{kavousi2018stochastic} \cite{fang2019optimum} \cite{fang2019adaptive} \\

                 &                        &  \cite{wang2020resilience} \cite{wu2021risk} \cite{alkhaleel2022risk} \cite{yang2018towards} \cite{ye2015resilience}  \cite{liu2021hierarchical} \\
     
                               \hdashline
\multirow{3}{*}{Miscellaneous}                     &Distance & \cite{akbari2017multi} \cite{iloglu2018integrated}\\

                 &Societal  & \cite{yang2021real}   \cite{aslani2022learn}     \\

     &Environmental  &   \cite{aslani2022learn}     \\

\hline
\end{tabular}
\label{obj}
\end{sidewaystable}

\subsection{Constraints}
Mathematical models attempt to capture technical limitations and rules (e.g., precedence) in CPS in terms of \textit{constraints}. The constraints for optimization models reviewed in this work are very diverse and problem-dependant. However, in this section, we provide a general scheme of a mathematical model for CPS to provide a practical guideline for researchers to build novel models based on the well-established models in the literature.

A general structure of an optimization model is presented in algorithm \ref{model}. Infrastructure networks provide a service (in terms of commodity flow) from supply nodes to their target customers (demand nodes). Therefore, the first set of constraints are \textit{flow balance constraints)} to preserve the balance between the inflow and outflow for demand, supply, and transshipment nodes (see \cite{sharkey2015interdependent} for a detailed explanation). The next constraints are the \textit{technical limitation constraints} to guarantee that the features such as capacity and upper limit for met demand of nodes are preserved during the optimization process (see \cite{lee2007restoration} for more information). These constraints can be in the form of equality and inequality equations based on the nature of technical limitations. As the problems in the literature considered isolated infrastructure systems and interdependent networks, these constraints can be different. Regardless of these differences, a set of \textit{dependency or Interdependency constraints} are included to capture the complex relationship among components (see \cite{cavdaroglu2013integrating} for detailed explanation). The \textit{resource allocation} and \textit{scheduling} set of constraints also aim to describe the assignment of repair crews and the sequence of recovery in the network, respectively (see \cite{nurre2012restoring} for more information). Finally, \textit{Application-focused constraints} are included to satisfy the specific requirements for the application area. One good example is the work of \cite{abessi2020new}, where they defined a set of specific equations to capture the assumption in power distribution systems that there are no switches on some lines in the system.
 







\begin{algorithm}[htbp]
\begin{align*}
&\begin{rcases}
    \sum x_{ijt}^{m} - \sum x_{jit}^{m}=  s_i^m \\
    \sum x_{ijt}^{m} - \sum x_{jit}^{m}=  0 \\
   \sum x_{ijt}^{m} - \sum x_{jit}^{m}= -v_{it}^m \\
   \vdots
\end{rcases}
\text{Flow balance constraints}\\
&\begin{rcases}
0\leq x_{ijt}^m \leq u_{ij}^{m} \\
 \sum v_{it}^m = d_i^{m}\\
\vdots
\end{rcases}
\text{Technical limitation constraints}\\
&\begin{rcases}
  0\leq d_{i}^m -v_{it}^m\leq (1-y_{n,j,t}^{m,i}) (d_{i}^m) \\
\sum x_{jht}^n \leq s_{j}^n y_{n,j,t}^{m,i}\\
\vdots\\
\end{rcases}
\text{Dependency or Interdependency constraints}\\
&\begin{rcases}
\sum \alpha_{kijs}^m \leq 1  \\
\vdots\\
\end{rcases}
\text{Resource allocation constraints}\\
&\begin{rcases}
 \beta_{ijt}^{m}-\beta_{ij(t-1)}^{m}=\sum \alpha_{kijt}^m \\
\vdots\\
\end{rcases}
\text{Scheduling constraints}\\
&\begin{rcases}
\sum x_{ij} \leq U\\
\sum y_{ij} = E\\
\vdots\\
\end{rcases}
\text{Application-focused constraints}\\
&\begin{rcases}
 v_{it}>0 \\
x_{jht} \in \{0,1\}\\
\vdots\\
\end{rcases}
\text{Decision variables}
\end{align*}
	\floatname{algorithm}{}
	\caption{\textbf {. General structure of a mathematical model for cyber-physical infrastructure system}} \label{model}
\end{algorithm}

\subsection{Network Size} \label{nsize}
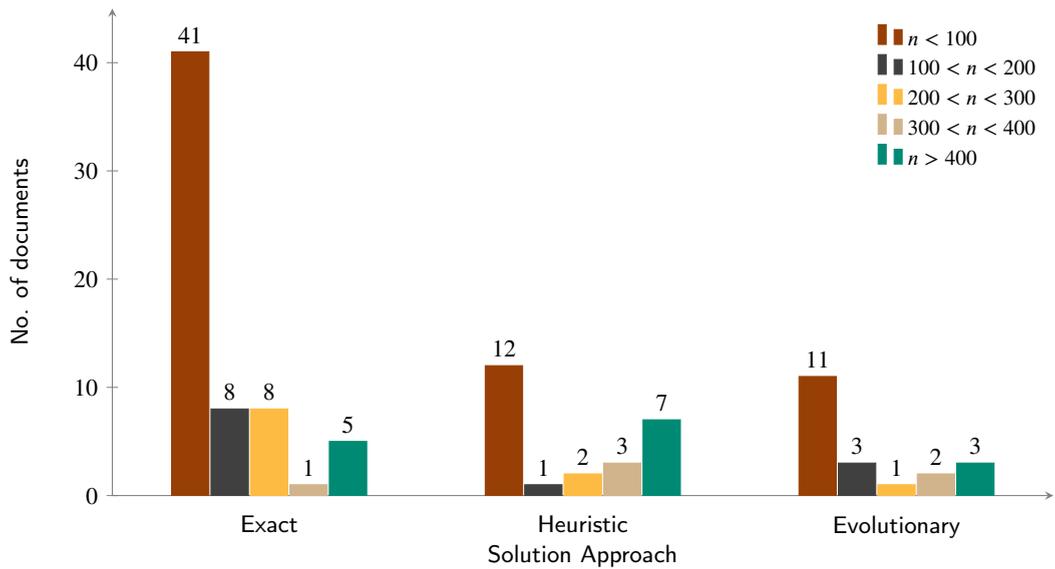
\begin{figure}[htpb]
\centering
\begin{tikzpicture}
    \begin{axis}[
        grid=none,
    axis lines=left,
        width  = 0.85*\textwidth,
        height = 8cm,
        major x tick style = transparent,
        ybar=2*\pgflinewidth,
        bar width=14pt,
        ylabel = {No. of documents},
        	y label style={at={(0,0.5)}},
         xlabel = {Solution Approach},
        symbolic x coords={Exact,Heuristic,Evolutionary},
        xtick = data,
        scaled y ticks = false,
        enlarge x limits=0.25,
        ymin=0,ymax=45,
        legend cell align=left,
        nodes near coords,
        every node near coord/.append style={text=black},
         legend style={draw=none,at={(0.9,0.98)},anchor=north, font=\footnotesize}
    ]
        \addplot[style={RawSienna,fill=RawSienna,mark=none}]
            coordinates {(Exact, 41) (Heuristic,12) (Evolutionary,11)};

        \addplot[style={darkgray,fill=darkgray,mark=none}]
             coordinates {(Exact,8) (Heuristic,1) (Evolutionary,3)};

        \addplot[style={Dandelion,fill=Dandelion,mark=none}]
             coordinates {(Exact,8) (Heuristic,2) (Evolutionary,1)};

        \addplot[style={Tan,fill=Tan,mark=none}]
             coordinates {(Exact,1) (Heuristic,3) (Evolutionary,2)};
             
        \addplot[style={PineGreen,fill=PineGreen,mark=none}]
             coordinates {(Exact,5) (Heuristic,7) (Evolutionary,3)};
        \legend{$n<100$,$100<n<200$,$200<n<300$,$300<n<400$,$n>400$}
        
    \end{axis}
\end{tikzpicture}
    \caption{No. of papers in different network size categories separated based on the solution methodology}
     \label{size}
\end{figure}

As the network expands to a city-scale level (a large number of nodes and arcs), the number of decision variables for mathematical models increases, and the consequent computational time rises exponentially \citep{cavdaroglu2013integrating}. Therefore, developing scalable and reliable solution algorithms for large-scale infrastructure networks is a demanding task from a purely technical perspective. As depicted in Figure \ref{size}, most developed optimization frameworks are designed for partial networks belonging to a much more extensive network. For example, \cite{sharkey2015interdependent}, \cite{baidya2017effective}, and \cite{kalinowski2015incremental} considered small and medium-sized networks (based on the definitions presented in section \ref{search}). In contrast, regarding the definition of a large-scale network, only a small share of studies including \cite{kasaei2016arc} and \cite{vodak2018modified} with approximately 1000, \cite{nurre2012restoring} with 760, and \cite{mohebbi2021decentralized} and \cite{aslani2022learn} with a city-scale interdependent network with around 12000 components developed heuristic and evolutionary approaches to tackle the computational burden of deterministic optimization problems. For stochastic optimization, \cite{gomez2019optimization} with 8400, and \cite{ccelik2015post} with 600 components are among the few studies focusing on a large-scale CPS network. 
   
\subsection{Failure Types} \label{failure}
A \textit{failure} can be broadly defined as an interruption of a CPS normal service level in providing the basic needs, such as communications, health, mobility, power, water, and sewer. The failure can be an immediate impact of an external factor such as cyber-attacks or aging of physical components. However, cascading failures are also common in these systems where a disruption in one CPS propagates into an interdependent CPS. Based on the source of failures, we categorized the failures in the literature into three main groups: 

\begin{itemize}
    \item \textbf{Functional Failures}: This type of failure is defined as the inability of a physical component to fulfill the intended function(s) to a standard performance level \citep{lee2007restoration}. The functional failures, originating from the internal mechanisms of cyber-physical systems, have been considered as random failures due to aging in several works such as \cite{bhuiyan2020stochastic} and \cite{karakoc2019community} and the cascading failure in others such as \cite{rahimi2022predictive} and \cite{alkhaleel2022risk}.  
    
    \item \textbf{Intentional Failures }: This type of failure is formally defined as a targeted failure by an intentional threat attempting to alter the normal operational level of a system \citep{fang2017optimizing}. Intentional failures have been investigated in the form of intentional attacks in several papers including \cite{ouyang2017mathematical} and \cite{fang2017optimizing}. The cyber-attacks, a special case of intentional failures that aim to interfere with cyber components, are also explored in some works such as \cite{wei2016quantifying} and \cite{moazeni2021formulating}. Finally, the localized attacks, spatially concentrated attacks, are investigated in \cite{ouyang2017spatial} and \cite{fang2019resilient}.
    
        \item \textbf{Natural Disaster}: This failure class is defined as any disruptions resulted from the immediate extreme event or the associated consequent events \citep{akbari2017multi}. The failures initiated by natural disasters are diverse in the literature with different triggering incidents such as earthquake in \cite{hu2013post}, flooding in \cite{sharkey2015interdependent}, Hurricane in \cite{baidya2017effective}, and windstorm in \cite{fang2019adaptive}. The failures resulting from natural disasters can be the immediate consequence of the incident or the secondary aftermath (e.g., the subsequent flooding of a hurricane). 

\end{itemize}

\subsection{Application Area} \label{application}
\begin{figure}[htbp]
	\centering
	\includegraphics[scale=0.7]{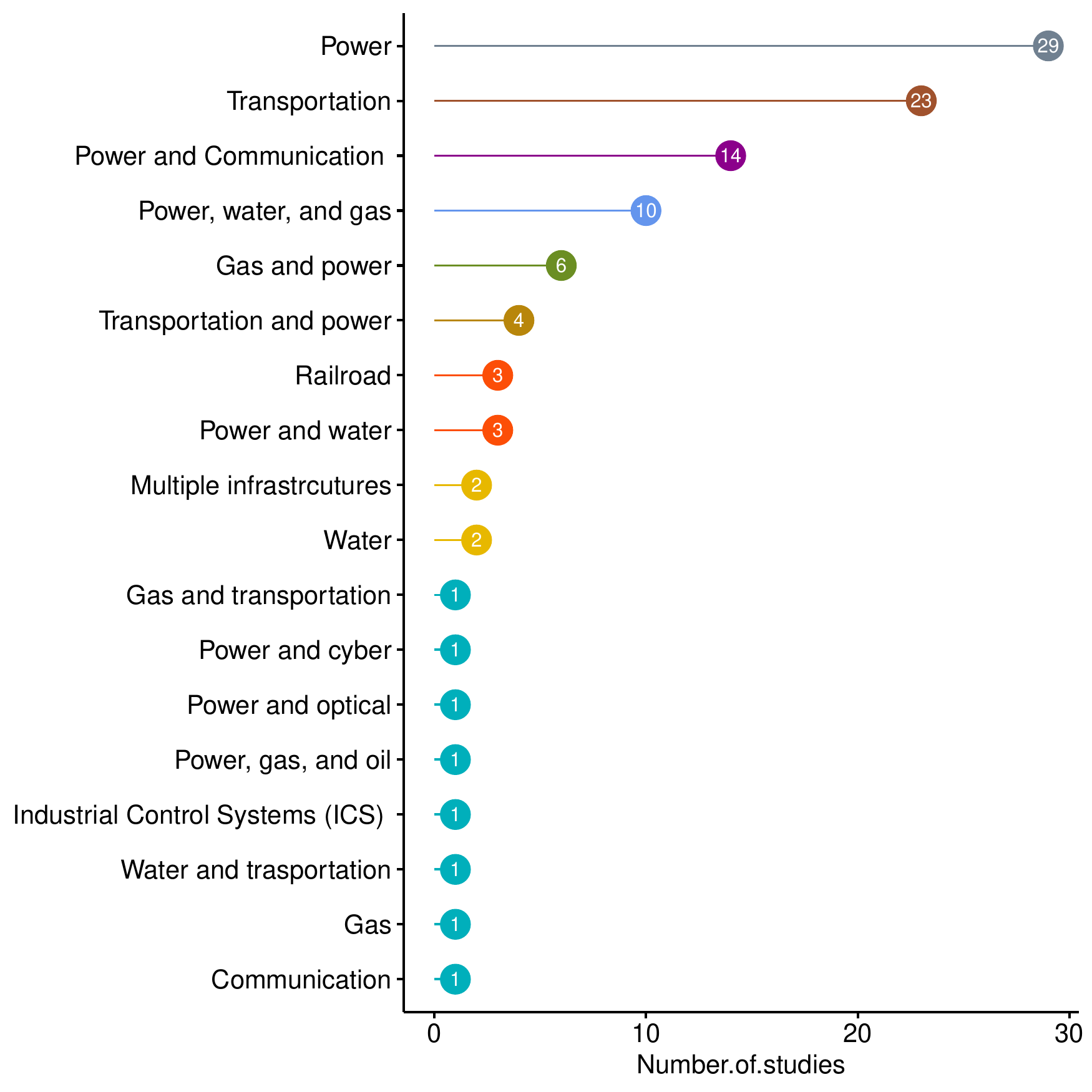}
	\centering
	\caption{Classification of application areas for optimization in the context of cyber-physical systems}
	\label{app}
\end{figure}

Fig. \ref{app} shows the type of cyber-physical infrastructure network considered for the optimization frameworks in the selected studies. The analysis shows that power and transportation networks are the most popular cyber-physical networks in the literature. This popularity can be stemmed from the availability of public data for these networks and their visibility. The numbers in this figure also reflect that the interdependent networks based on power and transportation and other sectors such as gas are also present numerously in the literature. On the other hand, gas, communication, and cyber networks (in the general sense) are among the least investigated application areas for optimization solution methods. In Fig. \ref{app} , the term multiple infrastructures in is refereed to studies that used more than 3 networks, like \textit{power, telecommunications, transportation, and wastewater} in \cite{sharkey2015interdependent} in their case studies. It should be noted that a share of selected documents, including \cite{kalinowski2015incremental} and \cite{duque2016network}, generated random networks for their case studies.      

`

\subsection{Interdependency Types} \label{interd}
The term \textit{interdependency} has been defined as bidirectional relationships between infrastructure systems where the state of one infrastructure affects the state of another network. Conceptually, the interdependencies are the complex connections among individual systems (networks) in a general system of systems functioning in an adaptive environment \citep{rinaldi2001identifying}. The interconnection in the context of CPS can emerge in four types; a) functional, b) cyber, c) geospatial (co-location), and d) logical \citep{mohebbi2020cyber}. The concept of functional interdependency refers to an interconnected network in which a component of one infrastructure requires services provided by another system to function correctly \citep{sharkey2015interdependent}. Co-location interdependency mirrors two infrastructure networks that can be affected by the same local disruptive event, where the initial failure propagates to the interdependent network and generates a cascading failure pattern \citep{rong2018optimum}. Cyber relationships among infrastructures reflect the situation where the state of one infrastructure system depends on information transmitted through the information infrastructure. Finally, logical interdependency also captures the relationship among the states of infrastructure systems by a mechanism that cannot be attributed to physical, cyber, or geospatial class \citep{ouyang2014review}.

\begin{figure}[htbp]
	\centering
	\includegraphics[scale=0.55]{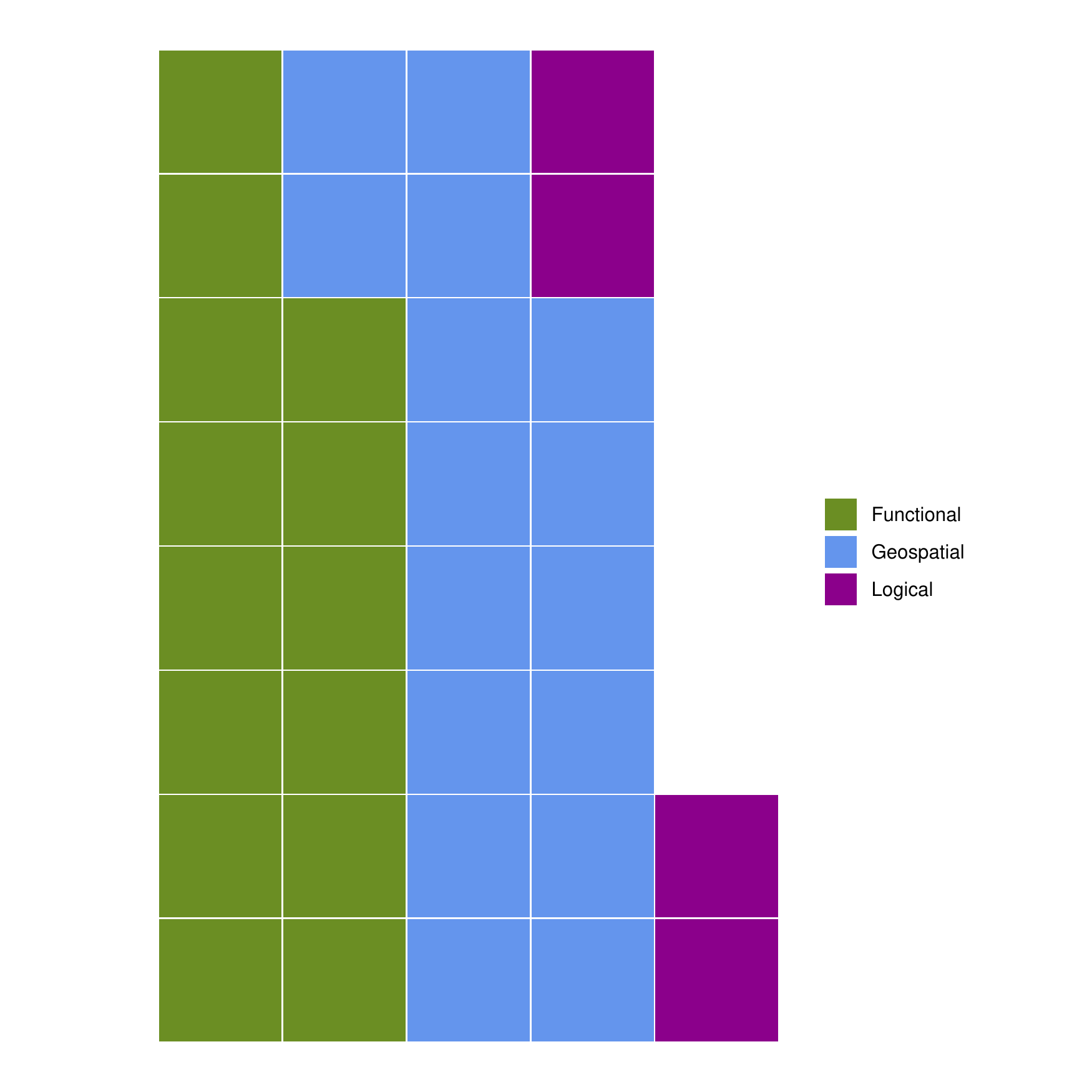}
	\centering
	\caption{Distribution (No. of studies) of interdependencies among infrastructure systems in the literature of optimization for cyber-physical systems}
	\label{inter}
\end{figure}

Fig. \ref{inter} shows the distribution of different types of interdependencies incorporated in the modeling for the selected studies. The functional and geospatial types are the most common as well-established methods (e.g., precedence constraints and spatial analysis) can capture these classes of relationships. On the other hand, only a few studies such as \cite{holden2013network}  and \cite{wang2020resilience} included the logical type considering the difficulty and vagueness in the definition of this category. Surprisingly, cyber interdependencies are missing in the literature of restoration of infrastructure systems. This absence can be attributed to the lack of a unanimously clear definition for cyber components and cyber relationships among physical elements of infrastructure networks.

\section{Discussion and Future Research} \label{discussion}

In this review, we attempted to provide a detailed analysis of the literature on optimization methods for CPS systems to answer two overarching questions: 1) \textit{What are the fundamental constructs of the developed mathematical models?}, and 2) \textit{What are the main trends for optimization-based solution algorithms for CPS?}. To answer the first question, we scrutinized the mathematical models in terms of the objective function, decision variables, and constraints. We critically analyzed the literature in each sub-category to provide concise yet insightful information for researchers and practitioners in developing customized models based on their specific needs. We also summarized the diverse pool of models in the literature to present an encapsulated mathematical model. Moving to the second question, we divided the selected documents in the review based on the modeling type (e.g., deterministic and stochastic) and solution method (e.g., exact, heuristic, and evolutionary). In each class, we identified the trends in the literature and the evolution over the time horizon between 2007 and 2022. In addition, we performed a bibliographic analysis to show how the methodology and application area are interconnected regarding keywords and publication outlets. The following sections, focusing on the significant gaps in the literature, prescribe possible future directions in the modeling and solution design perspectives.    

 \subsection{Scalable Solution Algorithms}
The term \textit{scalability} for optimization methods can be attributed to the ability of a solution algorithm to preserve an acceptable level of solution quality for large-scale instances of an optimization problem. The quality of solutions for an optimization method can be evaluated based on the optimality gap and computational time. While, generally speaking, the optimality gap for linear and convex models is zero (the solution is optimal), there are computational challenges for nonlinear and non-convex models. The optimality of a solution is not guaranteed for a nonlinear model when the non-convexity is present, implying there is a gap (usually presented in percentage format). In this context, the higher quality is attributed to a solution with a lower relative gap. Computational time is another feature to measure the solution quality, meaning that an optimization algorithm's ability to reach a good solution (not necessarily optimal) in lower computational time reflects the better quality. However, the two criteria should be assessed simultaneously to have an informative evaluation. In other words, a method's ability to provide solutions with a lower optimality gap in lower computational time is ideal. As a result, if a method can maintain this performance level for large-scale problem instances, we can label the method as a scalable solution algorithm. For example, \citep{xavier2021learning} developed such scalable hybrid ML-optimization method for unit commitment problems in power networks, showing that the large-scale instances can be solved on average 4.3 times faster with the optimality guarantee compared to conventional optimization methods. 

Specifically, the scalability of optimization methods is a major concern for decision-makers when dealing with city-scale CPS with several thousands of cyber-physical components. As the literature analysis revealed (section \ref{nsize}), the majority of developed optimization methods are designed for small (with lower than 100 overall components) and medium-scale (with lower than 400 overall components) networks. On the other hand, only a small share of studies approached large-scale networks and even fewer discussed the scalability feature. In specific, \citep{morshedlou2021heuristic} is the only paper that mentioned developing a local search heuristic to generate reliable solutions for large-scale problems. The current trends in the literature and the expansion of widespread CPS to modern communities reflect that developing scalable optimization methods to facilitate timely decision-making needs more attention and has remained an open search avenue.
    
\subsection{Hybrid Machine Learning and Optimization Algorithms}
A global optimization process includes exploring search space from an initial solution to find an optimal or near-optimal final solution. During this search mechanism, the optimization of the current problem is restarted from the beginning multiple times (using several heuristics to guide the search), and the algorithm collects some information about each iteration \citep{alvarez2017machine}. While these pieces of information contain valuable insights about the problem and recognizable patterns in the search space, the conventional optimization solution methods (e.g., using commercial solvers) usually discard these aspects of the search process \citep{lodi2020learning}. However, by identifying the sources of information, exploiting the data, and informing machine learning algorithms, the optimization algorithms can be guided intelligently to search the more promising areas and reach the final solution (optimal or near-optimal based on the nature of the algorithm) more efficiently. As a result, future research needs to focus on the interface of ML and optimization algorithms to exploit the values information and improve the scalability of solution methods.  

In the case of learning and optimization integration, all main learning streams (supervised, unsupervised, and reinforcement learning) can be embedded in the optimization algorithms for both deterministic and stochastic models. The choice of learning approach depends on the ability to collect relevant training data during the optimization procedure and the scope. If the optimization framework is known and the goal is to solve a class of problems (e.g., offline learning), the output of each problem instance can be used to guide a supervised/unsupervised learning process (i.e., training data) to solve the more challenging problems. A practical example is the work of \cite{xavier2021learning}, where they adopted two classical supervised learning methods for binary classification, k-nearest neighbors and support vector machines, to learn from the optimization output of small-scale problems and accelerate the solution process for large-scale cases of a deterministic MIP model. \cite{misra2021learning} took one further step by designing an algorithm to learn the relevant active constraints from training samples (e.g., input parameters and the corresponding optimal solution) without any restrictions on the problem type, problem structure, or probability distribution of the input parameters.

On the other hand, in the online learning framework, the learning and optimization happen simultaneously as solving a problem. In general, the choice of supervised/reinforcement learning for both deterministic and stochastic models in frameworks such as hyper-heuristic has been an active research direction in global optimization. For example, \cite{zhang2021deep} proposed a deep reinforcement learning-based hyper-heuristic framework to leverage a data-driven heuristic selection module on parameter-controlled low-level heuristics, to improve their handling of uncertainties while optimizing across combinatorial optimization problems.

There are other applications for ML-optimization integration for stochastic models (two-stage and multi-stage), where the uncertainty of parameters is presented in the form of scenarios. The large number of scenarios and computational challenges in Benders decomposition are the main bottlenecks of conventional methods. While supervised ML techniques can be implemented to cluster (or bundle) the scenarios, unsupervised ML algorithms can intelligently guide the decomposition algorithms and alleviate the computational burden. In the first research line,  \cite{jiang2021soft} presented a clustering method based on Fuzzy C-Means and Gaussian Mixture Models to bundle large-scale scenarios in a progressive hedging heuristic solution algorithm for multi-stage stochastic models. \cite{defourny2013scenario} also proposed a hybrid strategy based on ML and statistical inference from a given scenario-tree solution and out-of-sample simulation for multi-stage stochastic optimization problems. For the second direction, \cite{jia2021benders} developed a learning-enhanced Benders decomposition (augmented with a support vector machine classifier) to only generates promising cuts for solving two-stage stochastic programs with complete recourse based on finite samples of the uncertain parameters.

\subsection{Modular Optimization Framework}    
The literature of optimization solution algorithms for CPS includes a pool of various methods, including direct optimization (e.g.,  \cite{lee2007restoration}), heuristic search (e.g., \cite{wu2022allocation}), evolutionary methods (e.g., \cite{li2019joint}), and decomposition-based methods (e.g., \cite{alkhaleel2022risk}). Even though different solution algorithms have been developed in the literature for CPS, the methods are highly problem-dependant and cannot be generalized to other problem instances. In practice, decision-makers need to have easily reconfigurable solution methods to adapt to evolving situations and new problems with minimum technical efforts. In this regard, future research should investigate the modular design of the optimization framework, which makes the solution algorithms easily applicable to a wide range of problems and replaces the individual components based on the user's preferences.  

Object-oriented programming has emerged as a paradigm based on the concept of \textit{objects}, including data in the form of fields (often known as attributes or properties), and code, in the form of procedures (often known as methods). Inspired by this concept, one of the critical design requirements of optimization frameworks is the strict separation between algorithm and problem formulation. Hence, an optimization framework is divided into problem-specific (or domain) and general optimizer objects. While the former class of object should be tailored for each problem (in terms of solution representation and evaluation function), the latter group is independent of the domain. In addition, this object-oriented approach provides the flexibility of each object of the frameworks to be replaced with a similar one or compared in a computational way.

All in all, modular optimization frameworks are emerging general method-free solutions in other application areas. For example,  \cite{andersen2022mof} developed a modular optimization framework to create a flexible tool for users to easily incorporate new optimization algorithms, methods, or engineering design problems into the framework. In a similar work, \cite{kiss2020modular} designed modules that are replaceable and verifiable independently from each other to develop a framework capable of minimizing the life cycle environmental impact of a building design through automated optimization.

\subsection{State-of-the-art Optimization Methods}
As the summary of the literature in section \ref{class} revealed, the majority of solution algorithms for both deterministic and stochastic models are heavily biased towards conventional optimization methods. Even though we observe some efforts in works of \cite{mohebbi2021decentralized} and \cite{ghorbani2021decomposition} in developing novel hybrid strategies for exact solution algorithms, there is a need to shift to more novel methods. However, the flexibility of heuristic and evolutionary methods paves the way to leverage the data-driven approaches and develop innovative yet effective solution algorithms for deterministic and stochastic models. In specific, the promising results of \textit{hyper-heuristic} frameworks such as \cite{mosadegh2020stochastic} and the broad applications of learning in evolutionary algorithms such as \textit{cooperative co-evolution} frameworks (extensively reviewed in \cite{karimi2020learning}) should be considered fruitful avenues in developing future optimization frameworks.  

As another future direction, combining simulation techniques with heuristic/evolutionary algorithms known as \textit{simheuristics} aims to explore how problem-specific information can be used to enhance the solution method for combinatorial optimization problems with stochastic components (see \cite{juan2015review} for a comprehensive review). The stochastic elements can be present in the objective function or the constraints set. This novel integration promotes the use of risk-analysis criteria when evaluating alternative solutions. The simheuristics have been able to provide state-of-the-art solutions for combinatorial optimization problems in applications such ass vehicle routing, scheduling, manufacturing, system availability, and healthcare. Therefore, developing simheuristics for CPS systems is a promising direction to capture the uncertainty and provide practical solutions for decision-makers.

\subsection{Network Features for Optimization Methods}    
Networks are considered powerful tools to capture the complex dynamics of real-life systems. In the context of CPS, a spatial network is a graph in which the nodes are located in a space equipped with a specific metric, and the edges are spatial elements associated with nodes. The spatial features have been used mainly to capture the geospatial (co-location) interdependency among CPS networks in several studies such as \cite{rong2018optimum} and \cite{ge2019co}. The temporal feature of CPS, defined as the dynamic behavior of components over time, is another unique feature of networks. In this context, the temporal element of networks is utilized in \citep{rodriguez2022resilience} to capture the temporal interdependnecy. 

However, the spatiotemporal characteristics can be used at a higher level to inform optimization frameworks. For instance, the spatial information of cyber-physical components can be a decomposition criterion for focusing on decision variables decomposition. The temporal aspect also can be utilized to capture the uncertainty propagation in time-dependent optimization frameworks such as multi-stage stochastic models. In addition to spatiotemporal characteristics, \cite{evans2017optimization} argues that network measures (such as centrality) of individual elements (i.e., nodes or links) can quantify the associated connectivity profiles of elements and reflect how they are embedded in the network . Therefore, another promising area for optimizing CPS is exploiting such network features that can enrich the optimization algorithms designed specifically for large-scale networks.

\subsection{Formulation and Modeling Approaches}    

Despite the rich literature on mathematical modeling and conceptual frameworks related to CPS, there are certain deficiencies in methodology and application. As the first shortcoming, the optimization literature for CPS mainly focused on conventional modeling methods for capturing uncertainty, such as two-stage stochastic models. While this type of modeling is a popular choice in several application areas, there is a need to migrate towards more recent modeling approaches. In specific, Distributionally Robust Optimization (DRO) and Chance constrained optimization methods have been introduced as promising directions to capture uncertainty for stochastic models. For instance, \cite{zhao2020data} formulated the reliability of power-gas network in dealing with seismic attacks as a two-stage DRO model. The case studies showed that this novel modeling outperforms robust optimization and a single-stage optimization model to minimize the investment cost and expected economic loss. In addition, multi-stage stochastic modeling seems to be another modeling option that needs more attention, considering the nature of problems arising in CPS. While two-stage models are present in the works such as \cite{arab2015stochastic}, \cite{arif2018optimizing}, and \cite{abessi2020new}, the multi-stage stochastic models have not been explored in the literature deeply.

Other missing elements in the conceptual frameworks for CPS are the definition of cyber networks and cyber interdependency. The compound integration of cyber and physical levels for CPS shows the importance of identifying the cyber elements and cyber interdependencies to have more informative frameworks. While the analysis in section \ref{application} showed the term cyber network is used interchangeably with communication networks in several papers including \cite{huang2022optimization} and \cite{baidya2017effective}, the \cite{sandor2019cyber} defined industrial control systems as the cyber network. This variation in the definition of cyber networks shows the need to conceptually develop a clear understanding among researchers of what constructs a cyber network. Moreover, the output of section \ref{interd} revealed that the literature failed to capture and incorporate the cyber interdependency among CPS networks. \cite{mohebbi2020cyber} is among a few studies that attempted to identify cyber interdependencies for water and transportation infrastructure systems and how they may influence organizational resilience. They even recommended some design and operational strategies in this area. Nonetheless, these aspects are not reflected in conceptual frameworks and, consequently, optimization formulations in a comprehensive manner. 

The analysis of objective functions in section \ref{object} demonstrated that a large share of studies defined the objective function of the developed models as the \textit{resilience} of a CPS system. However, despite the dynamic nature of resilience in CPS, almost all studies in the literature simplified the term to a static measure (usually the proportional met demand). This conflict indicates that implementing dynamic modeling approaches for deterministic, multi-stage stochastic, and MDP models for stochastic formulation need further attention to capture the resilience evolution over time appropriately.    
    
Lastly, the summarized information in section \ref{failure} indicates that the literature has been mainly focused on standalone physical and cyber failures (under intentional failure class). However, the term CPS denotes that compound cyber-physical failures are major concerns in city-scale networks. Although some studies in the literature such as \cite{li2015bilevel}, \cite{he2021tri}, and \cite{jena2021design} focused on cyber-physical failures, future works should explore the impacts of such compound failures on CPS in a more detailed analysis. 

In summary, there are multiple future directions to extend the theory and application of mathematical modeling and optimization methods for CPS. As CPS become increasingly indispensable in modern communities, decision-makers seek to access practical optimization tools to prepare, respond, and adapt to evolving challenges.

\section*{Acknowledgements}
This work was supported by the Center for Resilient and Sustainable Communities (C-RASC) at George Mason University. Any opinions, findings, and conclusions or recommendations expressed in this material are those of the authors.

\section*{Disclosure statement}
The authors declare that they have no known competing financial interests or personal relationships that could have appeared to influence the work reported in this paper.

\bibliographystyle{apalike}
\bibliography{Review}

\end{document}